\documentclass[12pt]{article}
\usepackage{amssymb, amsmath, amsthm, amscd}
\usepackage[utf8]{inputenc}
\usepackage[all,cmtip]{xy}
\usepackage{enumitem}
\usepackage{setspace}
\usepackage{mathdots}

\usepackage[mathscr]{euscript}
\usepackage[colorlinks=true]{hyperref}
\hypersetup{colorlinks   = true}
\hypersetup{linkcolor=blue}
\usepackage{tikz}
\begin{document}

\mathchardef\mhyphen="2D
\newtheorem{The}{Theorem}[section]
\newtheorem{Lem}[The]{Lemma}
\newtheorem{Prop}[The]{Proposition}
\newtheorem{Cor}[The]{Corollary}
\newtheorem{Rem}[The]{Remark}
\newtheorem{Obs}[The]{Observation}
\newtheorem{SConj}[The]{Standard Conjecture}
\newtheorem{Titre}[The]{\!\!\!\! }
\newtheorem{Conj}[The]{Conjecture}
\newtheorem{Question}[The]{Question}
\newtheorem{Prob}[The]{Problem}
\newtheorem{Def}[The]{Definition}
\newtheorem{Not}[The]{Notation}
\newtheorem{Claim}[The]{Claim}
\newtheorem{Conc}[The]{Conclusion}
\newtheorem{Ex}[The]{Example}
\newtheorem{Fact}[The]{Fact}
\newtheorem{Formula}[The]{Formula}
\newtheorem{Formulae}[The]{Formulae}
\newtheorem{The-Def}[The]{Theorem and Definition}
\newtheorem{Prop-Def}[The]{Proposition and Definition}
\newtheorem{Lem-Def}[The]{Lemma and Definition}
\newtheorem{Cor-Def}[The]{Corollary and Definition}
\newtheorem{Conc-Def}[The]{Conclusion and Definition}
\newtheorem{Terminology}[The]{Note on terminology}
\newcommand{\C}{\mathbb{C}}
\newcommand{\R}{\mathbb{R}}
\newcommand{\N}{\mathbb{N}}
\newcommand{\Z}{\mathbb{Z}}
\newcommand{\Q}{\mathbb{Q}}
\newcommand{\Proj}{\mathbb{P}}
\newcommand{\Rc}{\mathcal{R}}
\newcommand{\Oc}{\mathcal{O}}
\newcommand{\Vc}{\mathcal{V}}
\newcommand{\Id}{\operatorname{Id}}
\newcommand{\pr}{\operatorname{pr}}
\newcommand{\rk}{\operatorname{rk}}
\newcommand{\del}{\partial}
\newcommand{\delbar}{\bar{\partial}}
\newcommand{\Cdot}{{\raisebox{-0.7ex}[0pt][0pt]{\scalebox{2.0}{$\cdot$}}}}
\newcommand\nilm{\Gamma\backslash G}
\newcommand\frg{{\mathfrak g}}
\newcommand{\fg}{\mathfrak g}
\newcommand{\Oh}{\mathcal{O}}
\newcommand{\Kur}{\operatorname{Kur}}
\newcommand\gc{\frg_\mathbb{C}}
\newcommand\hisashi[1]{{\textcolor{red}{#1}}}
\newcommand\dan[1]{{\textcolor{blue}{#1}}}

\begin{center}

{\Large\bf Partially Hyperbolic Compact Complex Manifolds}

\end{center}

\begin{center}

{\large Hisashi Kasuya and Dan Popovici}

\end{center}

\vspace{1ex}

\noindent{\small{\bf Abstract.} We propose and investigate two types, the latter with two variants, of notions of partial hyperbolicity accounting for several classes of compact complex manifolds behaving hyperbolically in certain directions, defined by a vector subbundle of the holomorphic tangent bundle, but not necessarily in the other directions. A key role is played by certain entire holomorphic maps, possibly from a higher-dimensional space, into the given manifold $X$. The dimension of the origin $\C^p$ of these maps is allowed to be arbitrary, unlike both the classical $1$-dimensional case of entire curves and the $1$-codimensional case introduced in previous work of the second-named author with S. Marouani. The higher-dimensional generality necessitates the imposition of certain growth conditions, very different from those in Nevanlinna theory and those in works by de Th\'elin, Burns and Sibony on Ahlfors currents, on the entire holomorphic maps $f:\C^p\longrightarrow X$. The way to finding these growth conditions is revealed by certain special, possibly non-K\"ahler, Hermitian metrics in the spirit of Gromov's K\"ahler hyperbolicity theory but in a higher-dimensional context. We then study several classes of examples, prove implications among our partial hyperbolicity notions, give a sufficient criterion for the existence of an Ahlfors current and a sufficient criterion for partial hyperbolicity in terms of the signs of two curvature-like objects introduced recently by the second-named author.}

\vspace{2ex}


\vspace{1ex}


\vspace{2ex}

\section{Introduction}\label{section:Introduction} In this paper, we continue the study of hyperbolic compact complex manifolds begun in [MP22a] and [MP22b] by relaxing the hyperbolicity requirement to some (i.e. not necessarily all) of the directions. This accommodates far larger classes of examples than in those references, while keeping the two-fold peculiarity of this point of view:

\vspace{1ex}

$(1)$\, possibly {\bf non-K\"ahler} manifolds are targeted, unlike the classical notions of hyperbolicity (due e.g. to Kobayashi, Brody, Gromov, etc) that only apply (at least conjecturally) to projective manifolds;

\vspace{1ex}

$(2)$\, holomorphic maps from possibly {\bf higher-dimensional} spaces (e.g. $\C^p$ with $p\geq 2$) are used in our study, unlike the classical entire curves defined on $\C$ that are involved in earlier notions of hyperbolicity.

\vspace{1ex}

In other words, the overall theme of this approach to hyperbolicity is the identification and the study of relations between the existence of special (possibly non-K\"ahler) Hermitian metrics and the constraints to which entire holomorphic maps from possibly higher dimensional spaces to compact complex manifolds are subject.

\vspace{2ex}

The starting point of this work is the standard fact that the universal covering space of Oeljeklaus-Toma (O-T) manifolds ([OT05]) is $\mathbb{H}^p\times\C^q$,  where $\mathbb{H}$ is the upper half-plane of $\C$, a typical hyperbolic manifold, while the complex Euclidean space $\C^q$ is as far as possible from being hyperbolic. The problem we set ourselves is to single out and study general properties of partial hyperbolicity that turn out to be displayed by quite a number of compact complex manifolds.

Let $X$ be a compact complex manifold with $\mbox{dim}_\C X=n\geq 2$. Suppose there exist $p\in\{1,\dots , n-1\}$ and a $C^\infty$ complex vector subbundle $E$ of rank $\geq p$ of the holomorphic tangent bundle $T^{1,\,0}X$. We will define two kinds of hyperbolicity properties that $X$ may have in the directions of $E$.

\vspace{1ex}

The former notion of partial hyperbolicity is of a metric nature and is reminiscent of Gromov's K\"ahler hyperbolicity of [Gro91], generalising it to $(p,\,p)$-forms on $X$ when $p$ may be larger than $1$. Specifically, suppose there exists a Hermitian metric $\omega=\omega_E + \omega_{nE}>0$, viewed as a positive definite $C^\infty$ $(1,\,1)$-form on $X$, that splits into the sum of two positive semi-definite $C^\infty$ $(1,\,1)$-forms such that $\omega_E$ is positive definite in the $E$-directions and $\omega_{nE}$ vanishes in the same directions. If the $(p,\,p)$-form $\Omega:=\omega_E^p/p!$ is $d$-closed and $\widetilde{d}$(bounded), we say that the manifold $X$ is {\bf partially $p$-K\"ahler hyperbolic} in the {\bf $E$-directions} (cf. Definition \ref{Def:partial_p-K-hyperbolic}). In the special case where $p=n-1$ and the $(n-1)$-K\"ahler hyperbolicity occurs in all the directions, we recover the earlier notion of {\it balanced hyperbolicity} introduced in [MP22a]. The terminology here is a nod to the notion of $p$-K\"ahler structure ($=$ a $d$-closed, $C^\infty$, strictly weakly positive $(p,\,p)$-form) introduced by Alessandrini and Bassanelli in [AB91]. However, the manifolds we study in this paper need not be $p$-K\"ahler, even when $p=1$, let alone K\"ahler. 

\vspace{1ex}

The latter notion of partial hyperbolicity for $X$ rules out the existence of a certain type of entire holomorphic maps $f:\C^p\longrightarrow X$ and thus reminds one of the Brody hyperbolicity of [Bro78] but, again, with $p$ possibly larger than $1$. The holomorphic maps $f$ that are excluded are non-degenerate at some point $x_0\in\C^p$, are {\it $E$-horizontal} (in the sense that the image of the differential map $df$ is contained in $E$) and, in a crucial departure from Brody's criterion of the case $p=1$, induce a relatively small growth for the volumes of the Euclidean balls in $\C^p$, measured against the (degenerate) metric pulled back from $X$ under $f$, as the radius tends to $\infty$. Moreover, we propose two variants of this latter notion of partial hyperbolicity according to whether this growth is subexponential (cf. Definition \ref{Def:subexp}), thus generalising the {\it divisorial hyperbolicity} of the earlier work [MP22a], or slightly faster (cf. the growth condition (\ref{eqn:Ahlfors-current})). We call such manifolds $X$ {\bf partially $p$-hyperbolic} in the {\bf $E$-directions} (cf. Definition \ref{Def:partial_p-hyperbolic}), respectively {\bf strongly partially $p$-hyperbolic} in the {\bf $E$-directions} (cf. Definition \ref{Def:strongly_partial-p-hyperbolic}).

Three of our results are summed up in the following statement (cf Theorems \ref{The:partial-p-hyperbolicity_implication}, \ref{The:partial-p-K-hyperbolicity-new-growth_implication} and \ref{The:strongly-partially-hyp_implies_partially-hyp}).

\begin{The}\label{The:introd_implications_partial} Let $X$ be a compact complex manifold with $\mbox{dim}_\C X=n$ and let $E\subset T^{1,\,0}X$ be a $C^\infty$ complex vector subbundle of rank $\geq p\in\{1,\dots , n-1\}$. The following implications hold: 

\vspace{2ex}  

\hspace{2ex} $(\star)$ $\begin{array}{lll} X \hspace{1ex} \mbox{is {\bf partially $p$-K\"ahler hyperbolic} in the {\bf $E$-directions}} &  &  \\
  \rotatebox{-90}{$\implies$} &  &  \\
  X \hspace{1ex} \mbox{is {\bf strongly partially $p$-hyperbolic} in the {\bf $E$-directions}} & & \\
  \rotatebox{-90}{$\implies$} &  &  \\
   X \hspace{1ex} \mbox{is {\bf partially $p$-hyperbolic} in the {\bf $E$-directions}.} &  &\end{array}$

\end{The}

To put these new definitions in the context of the existing ones, we spell out in the following diagram the various implications for the case where $p=n-1$.

\vspace{3ex}  

\hspace{2ex} $(\star\star)$ $\begin{array}{lll} X \hspace{1ex} \mbox{is {\bf K\"ahler hyperbolic}} & \implies & X \hspace{1ex} \mbox{is {\bf Kobayashi/Brody hyperbolic}} \\
 \rotatebox{-90}{$\implies$} &  & \rotatebox{-90}{$\implies$} \\
 X \hspace{1ex} \mbox{is {\bf balanced hyperbolic}} & \implies & X \hspace{1ex} \mbox{is {\bf divisorially hyperbolic}} \\
 \rotatebox{-90}{$\implies$} &  & \rotatebox{-90}{$\implies$}  \\
 X \hspace{1ex} \mbox{is {\bf partially $(n-1)$-K\"ahler hyperbolic}} & \implies & X \hspace{1ex} \mbox{is {\bf partially $(n-1)$-hyperbolic}}.\end{array}$

\vspace{3ex}

In $\S$\ref{section:examples_partial-hyp}, we exhibit four classes of compact non-K\"ahler  complex manifolds that have partial hyperbolicity properties: all Oeljeklaus-Toma manifolds [OT05] (cf. Proposition \ref{Prop:O-T_partially-s-K_hyp}), certain manifolds constructed very recently by Miebach and Oeljeklaus in [MO22] (cf. Proposition \ref{Prop:M-O}), a certain class of complex parallelisable solvmanifolds similar to those constructed in [Kas17] (cf. Proposition \ref{Prop:C-par-solvmanifolds}) and all compact Vaisman manifolds (cf. Proposition \ref{Prop:Vaisman}).

In $\S$\ref{section:Ahlfors}, we give a sufficient condition for a non-degenerate holomorphic map $f:\C^p\longrightarrow X$ to induce an {\it Ahlfors current}, an object that has played a key role in hyperbolicity issues since at least McQuillan's work [MQ98]. If $\omega$ is a Hermitian metric on $X$, one considers the bidegree-$(n-p, n-p)$-current $T_r$ on $X$ defined as the pushforward $f_\star[B_r]$ under $f$ of the current of integration on the Euclidean ball of radius $r>0$ centred at the origin of $\C^p$, normalised by a division by the volume of this ball with respect to the (possibly degenerate) metric $f^\star\omega$ on $\C^p$. When $r\to +\infty$, standard arguments enable one to extract a subsequence $(T_{r_\nu})_{\nu\in\N}$ converging in the weak topology of currents to a strongly positive, bidegree-$(n-p, n-p)$-current $T$ of unit mass w.r.t. $\omega$ on $X$. However, this current need not be {\it closed}. When it is, it is called an Ahlfors current. We prove that this is the case (cf. Theorem \ref{The:Ahlfors-current}) if $f$ satisfies the growth condition (\ref{eqn:Ahlfors-current}), the same that we then use in $\S$\ref{section:Ahlfors} to define our notion of {\it strong partial $p$-hyperbolicity}. Thus, an immediate consequence of Theorem \ref{The:Ahlfors-current} can be reworded as 

\begin{The}\label{The:introd_Ahlfors-current} Any compact complex manifold that is {\bf not strongly partially $p$-hyperbolic} in the directions of a given complex vector subbundle $E\subset T^{1,\,0}X$ of rank $\geq p$ carries an Ahlfors current.

\end{The}  

For this result, we drew inspiration from de Th\'elin's work [dT10]. However, our sufficient growth condition (\ref{eqn:Ahlfors-current}) is very different from his and is obtained by a different treatment, in the spirit of this paper, of certain integral estimates. This further enables us to relate it to our {\it subexponential growth condition} (cf. Proposition \ref{Prop:growth-notions_implication}).

It seems difficult to prove in concrete situations that a given compact complex manifold $X$ is either {\it strongly partially $p$-hyperbolic} or {\it partially $p$-hyperbolic} without proving beforehand that it has the stronger property of being {\it partially $p$-K\"ahler hyperbolic}. Indeed, metric structures are often easier to construct than entire holomorphic maps from some $\C^p$ are to rule out. This was one of our main motivations for introducing a metric hyperbolicity property that implies entire-map-based notions. This is also the reason why all our examples of partially hyperbolic manifolds exhibited in $\S$\ref{section:examples_partial-hyp} are {\it partially $p$-K\"ahler hyperbolic}. In particular, we do not have at this point in time any example of a manifold that is {\it partially $p$-hyperbolic}, but {\it non-strongly partially $p$-hyperbolic}. However, we hope such manifolds exist, possibly even among the submanifolds of certain projective manifolds or complex projective spaces $\C\Proj^n$, as the intuition suggests. Consequently, we believe both our entire-map-based notions of hyperbolicity will prove useful in future manifold classification considerations.  

\vspace{2ex}

In $\S$\ref{section:curvature}, we give a sufficient criterion for {\it strong partial $(n-1)$-hyperbolicity} in the directions of a complex vector subbundle $E\subset T^{1,\,0}X$ of co-rank $1$ on an $n$-dimensional compact complex manifold $X$ (cf. Theorem \ref{The:sufficient_metric-hyp}) in terms of the signs of two curvature-like objects: a function $f_\omega$ and an $(n-1,\,n-1)$-form $\star\rho_\omega$ uniquely associated with every Hermitian metric $\omega$ on $X$ via a construction introduced in [Pop23]. To this end, we define a notion of partial negativity (or {\it negativity in the $E$-directions}) for any $(n-1,\,n-1)$-form (cf. Definition \ref{Def:partial-negativity}) and observe that in several cases of explicit (and well-known) compact complex Hermitian manifolds $(X,\,\omega)$ the computations carried out in [Pop23] yield: \begin{eqnarray}\label{eqn:sign_curvature_assumption}f_\omega >0 \hspace{3ex} \mbox{and} \hspace{3ex}   \star\rho_\omega \hspace{1ex} \mbox{is {\bf negative in the $E$-directions}}.\end{eqnarray}

This is the hypothesis that we make to get our curvature-like criterion for partial hyperbolicity in the following reformulation of Theorem \ref{The:sufficient_metric-hyp}:

\begin{The}\label{The:introd_sufficient_metric-hyp} If there exists a Hermitian metric $\omega$ on $X$ satisfying property (\ref{eqn:sign_curvature_assumption}), $X$ is {\bf strongly partially $(n-1)$-hyperbolic} in the E-directions.

\end{The}

This is our analogue in the present context of a by now standard discussion. Indeed, the classical notions of hyperbolicity (e.g. in the sense of Kobayashi [Kob67] or Brody [Bro78]) are well known to be implied by various negativity assumptions on various curvature tensors, forms or functions associated with a given Hermitian metric on a complex manifold. For example, according to Theorem 3.8 in [Kob67], one has:

\vspace{1ex}

{\it A Hermitian manifold whose holomorphic sectional curvature is bounded above by a negative constant is Kobayashi hyperbolic.}

\vspace{1ex}

Similar classical results involving other types of curvature can be found in Kobayashi's book [Kob70, Chapter 2, $\S3$ and $\S4$]. Their proofs  make use of various forms of the maximum principle.

\vspace{1ex}

In the examples of partially hyperbolic manifolds that we exhibit in $\S$\ref{section:examples_partial-hyp} and where explicit computations can be performed, we analyse two types of curvature negativity/positivity that point to analogies with the classical setting:

\vspace{1ex}

-certain holomorphic sectional curvatures of a class of partially hyperbolic manifolds that includes the Oeljeklaus-Toma manifolds are found to be negative constants (cf. $\S$\ref{hol-sectional-curvature_O-T});

\vspace{1ex}

-the curvature form of the canonical bundle of each of the manifolds mentioned above is found to be at least semi-positive in all directions and positive definite in certain directions at every point (cf. $\S$\ref{curvature-K_X_O-T}).

\vspace{1ex}

Based on these special cases, one may wonder whether:

\vspace{1ex}

-by analogy with Kobayashi's $K_X$-ampleness conjecture for the standard Kobayashi/Brody hyperbolic manifolds $X$, some kind of (partial) (semi-)positivity of $K_X$ can be proved to hold for every {\it partially $p$-K\"ahler hyperbolic} compact complex manifold $X$;

\vspace{1ex}

-by analogy with classical results, the negativity of some classical curvatures implies one or more of our hyperbolicity properties introduced in this paper, despite the fact that in $\S$\ref{section:curvature} we felt the curvature-like objects introduced in [Pop23] to be best suited to the general case.

\vspace{2ex}

{\bf Further contextualisation of this paper.} We point out two examples of earlier works by other authors in order to emphasise the peculiarities of ours.

\vspace{1ex}

$\bullet$\, Demailly extensively used (see e.g. [Dem21]) the notion of {\it directed varieties} $(X,\,V)$ involving an irreducible closed analytic subspace $V\subset T^{1,\,0}X$ of the total space of the holomorphic tangent bundle of a complex manifold $X$ such that each fibre $V_x:=V\cap T^{1,\,0}_xX$ is a vector subspace of $T^{1,\,0}_xX$ and the map $X\ni x\mapsto\mbox{dim}_\C V_x$ is Zariski lower semicontinuous.

This is, of course, similar to our consideration of a complex vector subbundle $E$ of $T^{1,\,0}X$, but there are significant differences: Demailly's manifold $X$ is projective (while ours need not even be K\"ahler); Demailly works with entire curves $f:\C\longrightarrow X$ tangent to $V$ (while we deal with $E$-horizontal maps $f:\C^p\longrightarrow X$ by allowing $p$ to exceed $1$); no growth condition needs to be imposed on the entire curves used in Demailly's or the classical hyperbolicity theory thanks to Brody's reparametrisation lemma that only holds in complex dimension $1$.

In fact, a key aspect of our results is the identification of appropriate growth conditions that need to be placed on holomorphic maps $f:\C^p\longrightarrow X$ when $p\geq 2$ in order to sharpen our hyperbolicity notions and make them relevant. (As explained in [MP22a], if no growth conditions are imposed, the theory would miss many compact non-K\"ahler manifolds that enjoy hyperbolicity properties, so it would be highly unsatisfactory.) This is where the special Hermitian metrics come in: they led us to the realisation that, in replacing the K\"ahler metrics used by Gromov in his definition of the K\"ahler hyperbolicity by other types of metrics (e.g. the balanced metrics used in [MP22a] and [MP22b]), we would be able to spot the appropriate higher-dimensional substitute for Brody's lemma in the form of our growth conditions. These conditions involve volumes of balls in $\C^p$ and differ from the growth conditions used in the Nevanlinna theory that depend on the Nevanlinna characteristic $T_f(r)$.  

\vspace{1ex}

$\bullet$\, Burns and Sibony studied Ahlfors currents induced by holomorphic maps from a complex manifold $X$ of any dimension $k$ to a compact {\it K\"ahler} manifold $Y$ in [BS12]. However, the growth condition we impose on our maps in order to guarantee the existence of an Ahlfors current is different from both the one used in [dT10] and the one of [BS12]. It is more in tune with our purposes in this paper and enables us to find yet another notion of partial hyperbolicity.

In other words, much like the adaptation to our non-K\"ahler context of Gromov's notion of K\"ahler hyperbolicity, the use of Ahlfors currents under an assumption different from those of [dT10] and [BS12] becomes a tool, rather than a goal, for us.

\vspace{2ex}

\noindent {\it Acknowledgements}. Work on this paper began during the second-named author's visit to the Osaka University in January-February 2023 at the invitation of the first-named author. The former wishes to thank the latter and the entire Mathematics Department for hospitality and support. Both authors are grateful to S. Dinew for very useful comments on a preliminary version of this paper.

\section{Two notions of partial hyperbolicity}\label{section:partial-hyp_two-definitions}

In this section, we give the main definitions and prove a key implication between two of our notions. 

\subsection{Partially $p$-K\"ahler hyperbolic compact complex manifolds}\label{subsection:partially-p-K-hyperbolic}

Let $(X,\,\omega)$ be a complex $n$-dimensional manifold equipped with a Hermitian metric. We denote by $\pi_X:\widetilde{X}\longrightarrow X$ the universal covering map of $X$. Let $k\in\{0,\dots , 2n\}$ and let $\alpha$ be a $C^\infty$ differential form of degree $k$ on $X$. Recall that, according to [Gro91], $\alpha$ is said to be:

\vspace{1ex}

(i)\, $d(bounded)$ is there exists a $C^\infty$ form $\beta$ of degree $k-1$ on $X$ such that $\beta$ is bounded with respect to $\omega$ and $\alpha = d\beta$.

\vspace{1ex}

(ii)\, $\tilde{d}(bounded)$ if the lift $\pi_X^\star\alpha$ of $\alpha$ to the universal cover is $d(bounded)$ on $\widetilde{X}$ with respect to the metric $\widetilde\omega:=\pi_X^\star\omega$, the lift of $\omega$ to $\widetilde{X}$.

\vspace{2ex}

The first hyperbolicity notion that we propose in this paper for compact Hermitian manifolds is described in the following

\begin{Def}\label{Def:partial_p-K-hyperbolic} Let $X$ be a compact complex manifold with $\mbox{dim}_\C X =n\geq 2$. Let $p\in\{1,\dots , n-1\}$.

 If there exist:
  
\vspace{1ex}

\hspace{-5ex}(a)\, a $C^\infty$ complex vector subbundle $E\subset T^{1,\,0}X$ of rank $\geq p$ of the holomorphic tangent bundle of $X$; 

\vspace{1ex}

\hspace{-5ex}(b)\, positive semi-definite $C^\infty$ $(1,\,1)$-forms $\omega_E\geq 0$ and $\omega_{nE}\geq 0$ on $X$

\vspace{1ex}

\hspace{-5ex} with the following properties:

\vspace{1ex}

(i)\, the $C^\infty$ $(1,\,1)$-form $\omega:=\omega_E + \omega_{nE}$ is positive definite on $X$;

\vspace{1ex}

(ii)\, $\omega_E(x)(\xi,\,\bar\xi)>0$ for every point $x\in X$ and every $(1,\,0)$-tangent vector $\xi\in E_x$ lying in the fibre of $E$ over $x$;

\vspace{1ex}

(iii)\, $\omega_{nE}(x)(\xi,\,\bar\eta)=0$ for every point $x\in X$ and all $(1,\,0)$-tangent vectors $\xi,\,\eta\in E_x$ lying in the fibre of $E$ over $x$;

\vspace{1ex}

(iv)\, the $C^\infty$ $(p,\,p)$-form $\Omega:=\omega_E^p/p!$ is {\bf $d$-closed} and $\tilde{d}({\bf bounded})$ on $(X,\,\omega)$.

\vspace{1ex}

\noindent the manifold $X$ is said to be {\bf partially $p$-K\"ahler hyperbolic} in the {\bf $E$-directions}, the triple $(E,\,\Omega,\,\omega = \omega_E + \omega_{nE})$ is called a {\bf partially $p$-K\"ahler hyperbolic structure} on $X$ and $E$ is called the {\bf horizontal vector bundle}.

\end{Def}

The role of the subbundle $E$ in the above definition is to indicate the directions along which $\omega_E$ is positive definite. The hyperbolicity is required to hold only in these directions, accounting for its ``partial'' character. On the other hand, if $\omega_E$ happens to be positive definite in all the directions, one can choose $E = T^{1,\,0}X$ and $\omega_{nE}=0$. Then, the manifold $X$ is called {\bf $p$-K\"ahler hyperbolic}.

\begin{Obs}\label{Obs:foliation_p-K-hyperbolic} Let $X$ be a compact complex manifold with $\mbox{dim}_\C X =n\geq 2$ and let $p\in\{1,\dots , n-1\}$. Suppose there exists a {\bf Frobenius integrable} $C^\infty$ complex vector subbundle $E\subset T^{1,\,0}X$ of rank $> p$ such that $X$ is {\bf partially $p$-K\"ahler hyperbolic} in the {\bf $E$-directions}.

  Then, for any compact leaf $Y\subset X$ of the holomorphic foliation induced by $E$, the manifold $Y$ is {\bf $p$-K\"ahler hyperbolic}.

\end{Obs}

\noindent {\it Proof.} Using the notation in Definition \ref{Def:partial_p-K-hyperbolic}, we see that $T^{1,\,0}_yY = E_y$ for every $y\in Y$ implies that the restriction to $Y$ of $\omega_E$ is a $C^\infty$ positive definite $(1,\,1)$-form, hence induces a Hermitian metric, on the compact complex manifold $Y$. Moreover, the $C^\infty$ $(p,\,p)$-form $\Omega_{|Y}=(\omega_{E|Y})^p/p!$ is {\bf $d$-closed} and $\tilde{d}({\bf bounded})$ on $(Y,\,\omega_{|Y})$. The contention follows. \hfill $\Box$

\vspace{2ex}

  When $p\leq n-2$, the $d$-closedness of $\Omega$ implies that $d\omega_E = 0$ at every point $x\in X$ where $\omega_E(x)>0$. Indeed, if $x$ is such a point, $\omega_E$ is positive definite, hence defines a Hermitian metric, on a neighbourhood $U$ of $x$. We have $\omega_E^{p-1}\wedge d\omega_E = 0$ on $U$, hence $d\omega_E = 0$ since the pointwise multiplication map $\omega_E^{p-1}\wedge\cdot:\Lambda^3T^\star X\longrightarrow\Lambda^{2p+1}T^\star X$ is injective due to the fact that $p-1\leq n-3$ when $p\leq n-2$. (Recall that for any Hermitian metric $\gamma$ on an $n$-dimensional complex manifold $X$ and for any non-negative integer $k\leq n$, the pointwise multiplication map $\gamma^p\wedge\cdot:\Lambda^kT^\star X\longrightarrow\Lambda^{k+2p}T^\star X$ is injective when $p\leq n-k$ and surjective when $p\geq n-k$.)

  In particular, if $\omega=\omega_E> 0$ on $X$, then:

  \vspace{1ex}

  $\bullet$ $\omega$ is a K\"ahler metric on $X$ if $p\leq n-2$;

  \vspace{1ex}

  $\bullet$ $\omega$ is a {\it balanced hyperbolic} metric on $X$ in the sense of [MP22a] (recalled below) if $p=n-1$.

  \vspace{1ex}

 \noindent Note that a {\it partially $p$-K\"ahler hyperbolic} $n$-dimensional manifold $(X,\,E,\,\Omega)$ need not be K\"ahler.

\vspace{2ex}

In the case $p=n-1$, the notion of {\it partial $(n-1)$-K\"ahler hyperbolicity} introduced above generalises the notion of {\it balanced hyperbolicity} introduced in [MP22a]. The latter generalised, in turn, Gromov's notion of {\it K\"ahler hyperbolicity} introduced in [Gro91]. Recall that, according to [Gro91], a compact complex manifold $X$ is said to be {\it K\"ahler hyperbolic} if it carries a $\tilde{d}(bounded)$ K\"ahler metric $\omega$.

Meanwhile, according to [MP22a], a compact complex $n$-dimensional manifold $X$ is said to be {\it balanced hyperbolic} if it carries a balanced metric $\omega$ (namely, a $C^\infty$ positive definite $(1,\,1)$-form $\omega$ such that $d\omega^{n-1}=0$) such that $\omega^{n-1}$ is $\tilde{d}(bounded)$.

\subsection{Partially $p$-hyperbolic compact complex manifolds}\label{subsection:partially-p-hyperbolic}

 We will need a few preliminaries.

On the one hand, recall that if $X$ is an $n$-dimensional complex (not necessarily compact) manifold equipped with a Hermitian metric $\omega$, for a given integer $k\in\{0,\dots , n\}$, a $k$-form $v$ on $X$ is said to be {\it primitive} (w.r.t. $\omega$) if $\omega^{n-k+1}\wedge v =0$. This is known to be equivalent to $\Lambda_\omega v=0$, where $\Lambda_\omega$ is the adjoint of the Lefschetz operator $\omega\wedge\cdot$ w.r.t. the pointwise inner product $\langle\,\,,\,\,\rangle_\omega$ defined by $\omega$. In particular, all $k$-forms with $k\in\{0,\,1\}$ are primitive and so are all $(p,\,0)$-forms and all $(0,\,q)$-forms.

The following standard formula (cf. e.g. [Voi02, Proposition 6.29, p. 150]) for the Hodge star operator $\star = \star_\omega$ of $\omega$ acting on {\it primitive} forms $v$ of arbitrary bidegree $(p, \, q)$ will come in handy: \begin{eqnarray}\label{eqn:prim-form-star-formula-gen}\star\, v = (-1)^{k(k+1)/2}\, i^{p-q}\, \frac{\omega^{n-p-q}\wedge v}{(n-p-q)!}, \hspace{2ex} \mbox{where}\,\, k:=p+q.\end{eqnarray}

\vspace{1ex}

On the other hand, we will need preliminaries that are similar to those in [MP22a]. Let $f:\C^q\longrightarrow (X,\,\omega)$ be a holomorphic map to a compact complex Hermitian manifold with $n:=\mbox{dim}_\C X\geq 2$ and $1\leq q\leq n-1$. We will suppose that $f$ is non-degenerate at some point $x_0\in\C^q$, namely that $d_{x_0}f:\C^q\longrightarrow T_{f(x_0)}^{1,\,0}X$ is of maximal rank. Then, the set $$\Sigma_f:=\{x\in\C^q\,\mid\, f\,\,\mbox{is degenerate at}\,\,x\}$$ is a proper analytic subset of $\C^q$ and the $C^\infty$ $(1,\,1)$-form $f^\star\omega\geq 0$ on $\C^q$ is positive definite, hence it defines a Hermitian metric, on $\C^q\setminus\Sigma_f$. We will refer to $f^\star\omega$ as a {\it degenerate metric} on $\C^q$ with degeneration set $\Sigma_f$.

For any map as above and any $r>0$, we define the {\bf $(\omega,\,f)$-volume} of the ball $B_r\subset\C^q$ of radius $r$ centred at the origin to be \begin{eqnarray}\label{eqn:omega-f_vol_def}\mbox{Vol}_{\omega,\,f}(B_r):=\int\limits_{B_r}f^\star\omega_q,\end{eqnarray} where, for any $(1,\,1)$-form $\gamma\geq 0$ and any positive integer $p$, we set: $$\gamma_p:=\frac{\gamma^p}{p!}.$$

On the other hand, to define the areas of the spheres of $\C^q$ with respect to the degenerate metric $f^\star\omega$, we proceed as follows. For every $z\in\C^q$, let $\tau(z):=|z|^2$ be its squared Euclidean norm. At every point $z\in\C^q\setminus\Sigma_f$, the following equality holds by the definition of the Hodge star operator $\star_{f^\star\omega}$ induced by the metric $f^\star\omega$: \begin{eqnarray}\label{eqn:Hodge-star_degenerate}\frac{d\tau}{|d\tau|_{f^\star\omega}}\wedge\star_{f^\star\omega}\bigg(\frac{d\tau}{|d\tau|_{f^\star\omega}}\bigg) = f^\star\omega_q.\end{eqnarray} This implies that the $(2q-1)$-form $$d\sigma_{\omega,\,f}:=\star_{f^\star\omega}\bigg(\frac{d\tau}{|d\tau|_{f^\star\omega}}\bigg)$$ on $\C^q\setminus\Sigma_f$ is the {\it area measure} induced by $f^\star\omega$ on the spheres of $\C^q$. In other words, the restriction \begin{eqnarray}\label{eqn:area-measure_degenerate}d\sigma_{\omega,\,f,\,r}:=\bigg(\star_{f^\star\omega}\bigg(\frac{d\tau}{|d\tau|_{f^\star\omega}}\bigg)\bigg)_{|S_r}\end{eqnarray} is the area measure induced by the degenerate metric $f^\star\omega$ on the sphere $S_r\subset\C^q$ centred at the origin of radius $r$, for every $r>0$.

In particular, the {\bf $(\omega,\,f)$-area} of the sphere $S_r\subset\C^q$, namely the area w.r.t. $d\sigma_{\omega,\,f,\,r}$, is \begin{eqnarray}\label{eqn:omega-f_area_def}A_{\omega,\,f}(S_r):=\int\limits_{S_r}d\sigma_{\omega,\,f,\,r}>0, \hspace{5ex} r>0.\end{eqnarray}

We can now define the restriction that will be placed on the holomorphic maps into a given compact complex manifold $X$ to define another hyperbolicity property of $X$. This condition is the analogue of the one in [MP22a, Definition 2.3].

\begin{Def}\label{Def:subexp} Let $(X,\,\omega)$ be a compact complex Hermitian manifold with $\mbox{dim}_\C X =n\geq 2$. Let $q\in\{1,\dots , n-1\}$ and let $f:\C^q\longrightarrow X$ be a holomorphic map that is non-degenerate at some point $x_0\in\C^q$.

  We say that $f$ has {\bf subexponential growth} if the following two conditions are satisfied:

  \vspace{1ex}

  (i)\, there exist constants $C_1>0$ and $r_0>0$ such that \begin{eqnarray}\label{eqn:subexp_1}\int\limits_{S_t}|d\tau|_{f^\star\omega}\,d\sigma_{\omega,\,f,\,t}\leq C_1t\, \mbox{Vol}_{\omega,\,f}(B_t), \hspace{5ex} t>r_0;\end{eqnarray}

  (ii)\, for every constant $C>0$, we have: \begin{eqnarray}\label{eqn:subexp_2}\limsup\limits_{b\to +\infty}\bigg(\frac{b}{C} - \log F(b)\bigg) = +\infty,\end{eqnarray} where $$F(b):=\int\limits_0^b\mbox{Vol}_{\omega,\,f}(B_t)\,dt = \int\limits_0^b\bigg(\int\limits_{B_t}f^\star\omega_q\bigg)\,dt, \hspace{5ex} b>0.$$

\end{Def}

Since $X$ is compact, any two Hermitian metrics $\omega_1$ and $\omega_2$ on $X$ are comparable (in the sense that there exists a constant $A>0$ such that $(1/A)\,\omega_2\leq\omega_1\leq A\,\omega_2$). Therefore, the subexponential growth condition on the holomorphic maps into $X$ is independent of the choice of Hermitian metric on $X$.

\vspace{2ex}

 We will often be considering maps of the type described in the following

\begin{Def}\label{Def:horizontal-map} Let $X$ be an $n$-dimensional complex manifold. Suppose there exists a $C^\infty$ complex vector subbundle $E$ of $T^{1,\,0}X$ of rank $p\in\{1,\dots , n-1\}$.

  For any $q\in\{1,\dots , p\}$, a holomorphic map $f:\C^q\longrightarrow X$ is said to be {\bf $E$-horizontal} if, for every $x\in\C^q$, the image of its differential map $d_xf:\C^q\longrightarrow T_{f(x)}^{1,\,0}X$ at $x$ is contained in $E_{f(x)}$. 

\end{Def}

 We are now in a position to define our second notion of partial hyperbolicity.

\begin{Def}\label{Def:partial_p-hyperbolic} Let $X$ be a compact complex manifold with $\mbox{dim}_\C X =n\geq 2$. Let $p\in\{1,\dots , n-1\}$.

  If there exists a $C^\infty$ complex vector subbundle $E\subset T^{1,\,0}X$ of rank $\geq p$ of the holomorphic tangent bundle of $X$ such that there is no {\bf $E$-horizontal} holomorphic map $f:\C^p\longrightarrow X$ that is {\bf non-degenerate} at some point $x_0\in\C^p$ and has {\bf subexponential growth} in the sense of Definition \ref{Def:subexp}, the manifold $X$ is said to be {\bf partially $p$-hyperbolic} in the {\bf $E$-directions}.

\end{Def}

As in Definition \ref{Def:partial_p-K-hyperbolic}, the role of the subbundle $E$ is to indicate the directions along which the manifold $X$ is hyperbolic. On the other hand, if we can choose $E = T^{1,\,0}X$ in the above definition, $X$ is said to be {\bf $p$-hyperbolic}. If this last situation occurs for $p=n-1$, $X$ is {\it divisorially hyperbolic} in the sense of [MP22a]. The analogue of Observation \ref{Obs:foliation_p-K-hyperbolic} is the following

\begin{Obs}\label{Obs:foliation_p-hyperbolic} Let $X$ be a compact complex manifold with $\mbox{dim}_\C X =n\geq 2$ and let $p\in\{1,\dots , n-1\}$. Suppose there exists a {\bf Frobenius integrable} $C^\infty$ complex vector subbundle $E\subset T^{1,\,0}X$ of rank $> p$ such that $X$ is {\bf partially $p$-hyperbolic} in the {\bf $E$-directions}.

  Then, for any compact leaf $Y\subset X$ of the holomorphic foliation induced by $E$, the manifold $Y$ is {\bf $p$-hyperbolic}.

\end{Obs}

\noindent {\it Proof.} Suppose there exists a holomorphic map $f:\C^p\longrightarrow Y$ that is non-degenerate at some point $x_0\in\C^p$ and has subexponential growth in the sense of Definition \ref{Def:subexp}. Since $T^{1,\,0}_yY = E_y$ for every $y\in Y$, the $Y$-valued map $f$, when viewed as $X$-valued, is $E$-horizontal. But then, the existence of such a map $f$ contradicts the partial $p$-hyperbolicity assumption made on $X$ in the $E$-directions.  \hfill $\Box$

\subsection{Partial $p$-K\"ahler hyperbolicity implies partial $p$-hyperbolicity}\label{subsection:partial-p-hyperbolicity_implication} We now prove the last horizontal implication in diagram $(\star\star)$ of the introduction with $n-1$ replaced by an arbitrary $p\in\{1,\dots , n-1\}$.

\begin{The}\label{The:partial-p-hyperbolicity_implication} Let $X$ be an $n$-dimensional compact complex manifold and let $p\in\{1,\dots , n-1\}$.

  If there exists a $C^\infty$ complex vector subbundle $E\subset T^{1,\,0}X$ of rank $\geq p$ such that $X$ is {\bf partially $p$-K\"ahler hyperbolic} in the {\bf $E$-directions}, then $X$ is {\bf partially $p$-hyperbolic} in the {\bf $E$-directions}.

\end{The}

\noindent {\it Proof.} Let us suppose that a partial $p$-K\"ahler hyperbolic structure $(E,\,\Omega=(\omega_E)_p,\,\omega= \omega_E + \omega_{nE})$ exists on $X$. In particular, $X$ has a Hermitian metric $\omega = \omega_E + \omega_{nE}$ with the properties in Definition \ref{Def:partial_p-K-hyperbolic}. Let $\widetilde\Gamma$ be a $C^\infty$ form of degree $2p-1$ on $\widetilde{X}$ such that $\widetilde\Gamma$ is bounded w.r.t. $\widetilde\omega:=\pi_X^\star\omega$ and $$\pi_X^\star\Omega = \pi_X^\star(\omega_E)_p=d\widetilde\Gamma.$$

We will prove by contradiction that $X$ satisfies the conditions of Definition \ref{Def:partial_p-hyperbolic} w.r.t. the vector subbundle $E$ fixed above. Suppose there exists an $E$-horizontal holomorphic map $f:\C^p\longrightarrow X$ that is non-degenerate at some point $x_0\in\C^p$ and has subexponential growth in the sense of Definition \ref{Def:subexp}. Let $\Sigma_f\subset\C^p$ be the degeneracy set of $f$. Since $\C^p$ is simply connected, $f$ lifts to $\widetilde{X}$, so there exists a holomorphic map $\widetilde{f}:\C^p\longrightarrow\widetilde{X}$ such that $f=\pi_X\circ\widetilde{f}$. In particular, \begin{eqnarray*}f^\star\Omega = \widetilde{f}^\star(\pi_X^\star\Omega) = d(\widetilde{f}^\star\widetilde\Gamma)  \hspace{5ex} \mbox{on}\hspace{1ex} \C^p.\end{eqnarray*}

On the other hand, since $f$ is $E$-horizontal and $\omega = \omega_E + \omega_{nE}$ with $\omega_{nE}(x)(\xi,\,\bar\eta) = 0$ for every $x\in X$ and all $(1,\,0)$-vectors $\xi,\eta\in E_x$, we have: $(f^\star\omega_{nE})(u,\,\bar{v}) = \omega_{nE}(f_\star u,\,f_\star \bar{v}) = 0$ for all $(1,\,0)$-vector fields $u,v$ in $\C^p$. Hence \begin{eqnarray}\label{eqn:f-star_omega_E}f^\star\omega = f^\star\omega_E  \hspace{5ex} \mbox{on}\hspace{1ex} \C^p.\end{eqnarray}

We will need the following

\begin{Lem}\label{Lem:Gamma-tilde_bounded} The $(2p-1)$-form $\widetilde{f}^\star\widetilde\Gamma$ is bounded in $\C^p$ w.r.t. $f^\star\omega$.

\end{Lem}
 
\noindent {\it Proof.} For any tangent vectors $v_1, \dots , v_{2p-1}$ on $\C^p$, we have: \begin{eqnarray*}|(\widetilde{f}^\star\widetilde\Gamma)(v_1, \dots , v_{2p-1})|^2 = |\widetilde\Gamma(\widetilde{f}_\star v_1, \dots , \widetilde{f}_\star v_{2p-1})|^2 \leq C\,|\widetilde{f}_\star v_1|^2_{\widetilde\omega}\cdots|\widetilde{f}_\star v_{2p-1}|^2_{\widetilde\omega} = C\,|v_1|^2_{f^\star\omega}\cdots|v_{2p-1}|^2_{f^\star\omega},\end{eqnarray*} where the inequality follows, for some constant $C>0$ independent of the $v_j$'s, from the $\widetilde\omega$-boundedness of $\widetilde\Gamma$ and the last equality follows from $\widetilde{f}^\star\widetilde\omega = (\pi_X\circ\widetilde{f})^\star\omega = f^\star\omega$.  \hfill $\Box$

\vspace{2ex}

\noindent {\it End of proof of Theorem \ref{The:partial-p-hyperbolicity_implication}.} We will compute and estimate $\mbox{Vol}_{\omega,\,f}(B_r)$ in two ways, where $B_r\subset\C^p$ is the open ball of radius $r>0$ centred at the origin.

\vspace{1ex}

$\bullet$ On the one hand, applying the classical Fubini Theorem, we get: \begin{eqnarray}\label{eqn:Volume_Fubini}\mbox{Vol}_{\omega,\,f}(B_r) = \int\limits_{B_r}(f^\star\omega_p) = \int\limits_0^r\bigg(\int\limits_{S_t}d\mu_{\omega,\,f,\,t}\bigg)\,dt = \int\limits_{B_r}d\mu_{\omega,\,f,\,t}\wedge\frac{d\tau}{2t}, \hspace{5ex} r>0,\end{eqnarray} where we used the equality $\tau = t^2$ on the sphere $S_t\subset\C^p$ of radius $t$, which implies $d\tau = 2tdt$ on $S_t$, and where $d\mu_{\omega,\,f,\,t}$ is the positive measure on $S_t$ defined by the requirement: \begin{eqnarray*}\frac{1}{2t}\,d\mu_{\omega,\,f,\,t}\wedge(d\tau)_{|S_t} = (f^\star\omega_p)_{|S_t},  \hspace{5ex} t>0.\end{eqnarray*} 

Comparing this with (\ref{eqn:Hodge-star_degenerate}) (in which we take $q=p$), we infer the following equality on $\C^p\setminus\Sigma_f$: \begin{eqnarray}\label{eqn:comparison_measures}\frac{1}{2t}\,d\mu_{\omega,\,f,\,t} = \frac{1}{|d\tau|_{f^\star\omega}}\,d\sigma_{\omega,\,f,\,t}, \hspace{5ex} t>0.\end{eqnarray}

On the other hand, H\"older's inequality yields: \begin{eqnarray*}A_{\omega,\,f}^2(S_t) = \bigg(\int\limits_{S_t}d\sigma_{\omega,\,f,\,t}\bigg)^2\leq\bigg(\int\limits_{S_t}\frac{1}{|d\tau|_{f^\star\omega}}\,d\sigma_{\omega,\,f,\,t}\bigg)\,\bigg(\int\limits_{S_t}|d\tau|_{f^\star\omega}\,d\sigma_{\omega,\,f,\,t}\bigg), \hspace{5ex} t>0.\end{eqnarray*} 

Together with (\ref{eqn:Volume_Fubini}) and (\ref{eqn:comparison_measures}), this yields: \begin{eqnarray}\label{eqn:Volume_1}\nonumber\mbox{Vol}_{\omega,\,f}(B_r) & = & \int\limits_0^r\bigg(\int\limits_{S_t}\frac{1}{2t}\,d\mu_{\omega,\,f,\,t}\bigg)\,d\tau = \int\limits_0^r\bigg(\int\limits_{S_t}\frac{1}{|d\tau|_{f^\star\omega}}\,d\sigma_{\omega,\,f,\,t}\bigg)\,d\tau \\
  & \geq & 2\,\int\limits_0^r\frac{A_{\omega,\,f}^2(S_t)}{\int\limits_{S_t}|d\tau|_{f^\star\omega}\,d\sigma_{\omega,\,f,\,t}}\,t\,dt, \hspace{5ex} r>0.\end{eqnarray}

\vspace{1ex}

We stress that (\ref{eqn:Volume_1}) holds for any holomorphic map $f:\C^p\longrightarrow (X,\,\omega)$ that is non-degenerate at some point $x_0\in\C^p$ and takes values in any complex Hermitian manifold. This map need not be either $E$-horizontal or of subexponential growth.

\vspace{1ex}

$\bullet$ On the other hand, applying the Stokes Theorem we get equality (d) below: \begin{eqnarray}\label{eqn:Volume_Stokes}\nonumber\mbox{Vol}_{\omega,\,f}(B_r) & = & \int\limits_{B_r}f^\star\omega_p \stackrel{(a)}{=} \int\limits_{B_r}f^\star(\omega_E)_p \stackrel{(b)}{=} \int\limits_{B_r}\widetilde{f}^\star(\pi_X^\star\Omega) \stackrel{(c)}{=}\int\limits_{B_r}d(\widetilde{f}^\star\widetilde\Gamma) \stackrel{(d)}{=} \int\limits_{S_r}\widetilde{f}^\star\widetilde\Gamma \\
  & \leq & C\,\int\limits_{S_r}d\sigma_{\omega,\,f,\,r} = C\,A_{\omega,\,f}(S_r), \hspace{5ex} r>0,\end{eqnarray} where (a) follows from $f$ being $E$-horizontal via (\ref{eqn:f-star_omega_E}), equalities (b) and (c) are consequences of the partial $p$-K\"ahler hyperbolicity assumption, while the inequality follows from Lemma \ref{Lem:Gamma-tilde_bounded}.

\vspace{1ex}

$\bullet$ Putting together (\ref{eqn:Volume_1}) and (\ref{eqn:Volume_Stokes}), we get the first inequality below: \begin{eqnarray}\label{eqn:Volume_final}\nonumber\mbox{Vol}_{\omega,\,f}(B_r) & \geq & \frac{2}{C^2}\,\int\limits_0^r \mbox{Vol}_{\omega,\,f}(B_t)\,\frac{t\,\mbox{Vol}_{\omega,\,f}(B_t)}{\int\limits_{S_t}|d\tau|_{f^\star\omega}\,d\sigma_{\omega,\,f,\,t}}\,dt \\
& \geq & \frac{2}{C_1C^2}\,\int\limits_{r_0}^r\mbox{Vol}_{\omega,\,f}(B_t)\,dt:=\frac{2}{C_1\,C^2}\,\widetilde{F}(r), \hspace{5ex} r>r_0,\end{eqnarray} where the second inequality follows from part (i) of the subexponential growth assumption on $f$ (cf. Definition \ref{Def:subexp}) for some constants $C_1, r_0>0$ and the last equality constitutes the definition of a function $\widetilde{F}:(r_0,\,+\infty)\longrightarrow(0,\,+\infty)$. 

Now, setting $C_2:=2/(C_1C^2)>0$, differentiating $\widetilde{F}$ and using (\ref{eqn:Volume_final}), we get: $$\widetilde{F}'(r) = \mbox{Vol}_{\omega,\,f}(B_r) \geq C_2\,\widetilde{F}(r), \hspace{5ex} r>r_0.$$ Hence $$\frac{d}{dt}\bigg(\log \widetilde{F}(t)\bigg)\geq C_2, \hspace{5ex} t>r_0.$$ 

Fixing arbitrary reals $a,b$ such that $r_0<a<b$ and integrating w.r.t. $t\in[a,\,b]$, we get: $$-\log \widetilde{F}(a) \geq -\log \widetilde{F}(b) + C_2\,(b-a).$$ Now, fixing $a>r_0$ and letting $b=b_j\to +\infty$ for a sequence of reals $b_j$ such that $-\log \widetilde{F}(b_j) + C_2\,b_j\to +\infty$ as $j\to +\infty$ (such a sequence exists when $\widetilde{F}$ is replaced by $F$ thanks to part (ii) of the subexponential growth assumption on $f$ -- see Definition \ref{Def:subexp} -- hence it also exists for $\widetilde{F}$ because $\widetilde{F}(r) = F(r) - F(r_0) \leq F(r)$ for every $r>r_0$), we get $\widetilde{F}(a) = 0$. Letting $a>r_0$ vary, this means that $\mbox{Vol}_{\omega,\,f}(B_t) = 0$ for every $t>r_0$. This implies $f^\star\omega_p = 0$ on $\C^p$, contraditing the non-degeneracy asssumption on $f$ and the property $\omega>0$ on $X$. \hfill $\Box$

\section{Examples of partially hyperbolic manifolds}\label{section:examples_partial-hyp}

We will describe several classes of examples, one in each of the ensuing subsections.

\subsection{Partial hyperbolicity and curvature of Oeljeklaus-Toma manifolds}\label{subsection:examples:O-T} 
 Given positive integers $s$ and $t$, we consider the Lie  group $G$ defined as the semi-direct product $G=\R^s\ltimes_\phi(\R^s\oplus\C^t)$ via the map $\phi:\R^s\longrightarrow\mbox{Aut}(\R^s\oplus\C^t)$, \begin{eqnarray}\label{eqn:semi-direct-prod_map_def_O-T}\phi(x)=\mbox{diag}\bigg(e^{x_1}, \dots , e^{x_s}, e^{\psi_1(x)}, \dots , e^{\psi_t(x)}\bigg),\end{eqnarray} where $x=(x_1,\dots , x_s)\in\R^s$ and $\psi_1,\dots ,\psi_t:\R^s\longrightarrow\C$ are functions of the shape $$\psi_j(x)=\sum\limits_{k=1}^sa_{jk}\,x_k,  \hspace{5ex} j=1,\dots , t,$$ with constant coefficients $a_{jk}\in\C$, while $\mbox{diag}$ stands for the diagonal matrix whose diagonal entries are those indicated.

Our objects of study in this $\S$\ref{subsection:examples:O-T} will be solvmanifolds $X=G/\Lambda$, namely quotients of the solvable Lie group $G$ by lattices $\Lambda\subset G$.

Denoting by $x=(x_1,\dots , x_s)$ the variable in the first copy of $\R^s$ in $G=\R^s\ltimes_\phi(\R^s\oplus\C^t)$, by $y=(y_1,\dots , y_s)$ the variable in the second copy of $\R^s$ and by $z=(z_1,\dots , z_t)$ the variable of $\C^t$, we see that the dual of the $(1,\,0)$-part of the Lie algebra of $G$ is generated as a $\C$-vector space as $$(\fg^{1,\,0})^\star = \bigg\langle dx_1+ie^{-x_1}\,dy_1,\dots , dx_s+ie^{-x_s}\,dy_s,\,e^{-\psi_1(x)}\,dz_1,\dots , e^{-\psi_t(x)}\,dz_t\bigg\rangle = \langle\alpha_1,\dots , \alpha_s,\,\beta_1,\dots , \beta_t\rangle,$$ where we have set: \begin{eqnarray}\label{eqn:O-T_alpha_j-beta_k}\nonumber\alpha_j & = & dx_j+ie^{-x_j}\,dy_j, \hspace{5ex} j=1,\dots , s, \\
          \beta_k & = & e^{-\psi_k}\,dz_k, \hspace{5ex} k=1,\dots , t.\end{eqnarray} The forms $\alpha_1, \dots , \alpha_s, \beta_1, \dots , \beta_t$ are $C^\infty$ $(1,\,0)$-forms on $G$ that induce an {\it invariant complex structure} on $G$. They also induce $C^\infty$ $(1,\,0)$-forms (which we denote by the same symbols, for the sake of convenience) on the solvmanifold $X=G/\Lambda$ and define a complex structure thereon, for any lattice $\Lambda\subset G$.

          The Lie group $G$ equipped with this complex structure is biholomorphic to $\mathbb{H}^s\times\C^t$ since the map $$\R\ltimes_\varphi\R\longrightarrow\mathbb{H}, \hspace{5ex} (x,\,y)\mapsto(y,\,e^x),$$ is biholomorphic, where the map $\varphi:\R\longrightarrow\mbox{Aut}(\R)$ is defined by  $$\varphi(x)(y)=e^x\,y, \hspace{5ex} y\in\R.$$ 

We note that the class of complex solvmanifolds of this type contains the compact complex  manifolds associated with algebraic number fields introduced in [OT05], the so-called Oeljeklaus-Toma (O-T) manifolds. This point of view was given in [Kas13] and was subsequently used in the study of metrical and cohomological properties of O-T manifolds (see [FKV15], [Kas21], [Kas23], [Oti22], [ADOS22]). 

Specifically, let $K$ be a finite extension field of $\Q$ of degree $s + 2t$ admitting real embeddings $\sigma_{1},\dots \sigma_{s}$ into $\C$ and complex embeddings $\sigma_{s+1},\dots, \sigma_{s+2t}$ into $\C$ that satisfy the conditions $\sigma_{s+i}=\bar \sigma_{s+i+t}$ for $1\le i\le t$. Meanwhile, let $U$ be a free subgroup of rank $s$ of the group of units in the ring ${\mathcal O}_{K}$ of algebraic integers satisfying certain conditions obtained by  Dirichlet's unit theorem. By choosing suitable functions $\psi_j:\R^s\longrightarrow\C$ as above, we can associate with the pair $(K,\,U)$ a lattice $\Lambda$ in $G$ that is isomorphic to $U\ltimes{\mathcal O}_{K}$. The resulting solvmanifold $X=G/\Lambda$ is an O-T manifold. Conversely, every O-T manifold is obtained in this way.

\subsubsection{Partial hyperbolicity}\label{subsubsection:partial-hyp_O-T}

Straightforward computations yield: $$d\psi_j = \sum\limits_{k=1}^sa_{jk}\,dx_k = \sum\limits_{k=1}^sa_{jk}\,\mbox{Re}\,\alpha_k = \frac{\psi_j(\alpha) + \psi_j(\overline\alpha)}{2},  \hspace{5ex} j=1,\dots , t,$$ where we put $\psi_j(\alpha):=\sum_{k=1}^sa_{jk}\,\alpha_k$ for every $j=1,\dots , t$. Since $\psi_j(\alpha)$ is a $(1,\,0)$-form and $\psi_j(\overline\alpha)$ is a $(0,\,1)$-form, this means that $\partial\psi_j = \psi_j(\alpha)/2$ and $\bar\partial\psi_j = \psi_j(\overline\alpha)/2$. 

On the other hand, the definition of $\alpha_j$ in (\ref{eqn:O-T_alpha_j-beta_k}) implies, through straightforward calculations, the following equalities: $d\alpha_j = -ie^{-x_j}\,dx_j\wedge dy_j$ and $(1/2)\,(\alpha_j\wedge\overline\alpha_j) = -ie^{-x_j}\,dx_j\wedge dy_j$. Since $\alpha_j\wedge\overline\alpha_j$ is of type $(1,\,1)$, we deduce that the $(2,\,0)$-part of $d\alpha_j$ must vanish, hence \begin{eqnarray*}\partial\alpha_j =  0 \hspace{2ex}\mbox{and}\hspace{2ex} \bar\partial\alpha_j =  d\alpha_j = \frac{\alpha_j\wedge\overline\alpha_j}{2} \end{eqnarray*} for every $j=1,\dots , s$. 

As for the $\beta_k$'s, we get: \begin{eqnarray*}\partial\beta_k & = & -e^{-\psi_k}\,\partial\psi_k\wedge dz_k = -\frac{\psi_k(\alpha)}{2}\wedge\beta_k  \\
\bar\partial\beta_k & = & -e^{-\psi_k}\,\bar\partial\psi_k\wedge dz_k = -\frac{\psi_k(\overline\alpha)}{2}\wedge\beta_k\end{eqnarray*} for all $k=1,\dots , t$.

The solvmanifold $X=G/\Lambda$ is a compact complex manifold of dimension $n=s+t$. Let $\{e_1,\dots , e_s, f_1,\dots , f_t\}$ be the $C^\infty$ frame of the holomorphic tangent bundle $T^{1,\,0}X$ of $X$ dual to the $C^\infty$ frame $\{\alpha_1,\dots , \alpha_s, \beta_1,\dots , \beta_t\}$ of the holomorphic cotangent bundle $\Lambda^{1,\,0}T^\star X$. (The vector bundles $T^{1,\,0}X$ and $\Lambda^{1,\,0}T^\star X$ are $C^\infty$-trivial, but not holomorphically trivial.) Let $E$ and $F$ be the $C^\infty$-trivial vector subbundles of $T^{1,\,0}X$ generated by $e_1,\dots , e_s$ and, respectively, by $f_1,\dots , f_t$. They also have natural, but non-trivial, holomorphic vector bundle structures.   

A natural Hermitian metric on $X$ is defined by the positive definite $C^\infty$ $(1,\,1)$-form \begin{eqnarray}\label{eqn:O-T_metric}\omega = i\alpha_1\wedge\overline\alpha_1 + \dots + i\alpha_s\wedge\overline\alpha_s + i\beta_1\wedge\overline\beta_1 + \dots + i\beta_t\wedge\overline\beta_t>0.\end{eqnarray} Thus, $\omega = \omega_E + \omega_F$, where $\omega_E:=i\alpha_1\wedge\overline\alpha_1 + \dots + i\alpha_s\wedge\overline\alpha_s\geq 0$ and $\omega_F:=i\beta_1\wedge\overline\beta_1 + \dots + i\beta_t\wedge\overline\beta_t\geq 0$ are $C^\infty$ positive semi-definite $(1,\,1)$-forms on $X$. By construction, $\omega_E(\xi,\,\bar\xi)>0$ for every $(1,\,0)$-tangent vector $\xi$ of $E$ and $\omega_E(\eta,\,\bar\nu) = 0$ for all $(1,\,0)$-tangent vectors $\eta,\,\nu$ of $F$. The analogous property is satisfied by $\omega_F$ with $E$ and $F$ permuted.

%Let us now consider the following $C^\infty$ $(s,\,s)$-form on $X$: \begin{eqnarray*}\Omega = \frac{\omega_E^s}{s!} = i\alpha_1\wedge\overline\alpha_1\wedge\cdots\wedge i\alpha_s\wedge\overline\alpha_s = 2id\alpha_1\wedge\cdots\wedge 2id\alpha_s = d\bigg(2i\alpha_1\wedge 2id\alpha_2\wedge\cdots\wedge 2id\alpha_s\bigg).\end{eqnarray*} 

Let us now consider the following $C^\infty$ $(1,\,1)$-form on $X$: \begin{eqnarray*}\Omega = \omega_E= d\bigg(2i\alpha_1+\cdots +2i\alpha_s\bigg).\end{eqnarray*} Note that $\Omega$ is not a Hermitian metric on $X$ since it is only $\geq 0$. Since $\Omega$ is $d$-exact on the {\it compact} complex manifold $X$, its lift to the universal cover $G\simeq\mathbb{H}^s\times\C^t$ is $d(bounded)$ with respect to the lift of the metric $\omega$. In particular, the triple $(E,\,\Omega,\,\omega= \omega_E + \omega_F)$ is a partially $1$-K\"ahler hyperbolic structure on $X$. 

We have thus proved the following

\begin{Prop}\label{Prop:O-T_partially-s-K_hyp} For every positive integers $s$ and $t$ and every simply connected solvable real Lie group $G=\R^s\ltimes_\phi(\R^s\oplus\C^t)$ defined as a semi-direct product via a map $\phi$ of the type (\ref{eqn:semi-direct-prod_map_def_O-T}), the solvmanifold $X=G/\Lambda$ obtained as the quotient of $G$ by any co-compact lattice $\Lambda$ is {\bf partially $1$-K\"ahler hyperbolic}. In particular every Oeljeklaus-Toma manifold is {\bf partially $1$-K\"ahler hyperbolic}.
\end{Prop}

\subsubsection{Holomorphic sectional curvatures}\label{hol-sectional-curvature_O-T}

We will now compute the holomorphic sectional curvatures of solvmanifolds in the $E$-directions (i.e. the hyperbolic directions) by continuing the above computations. Let $D=D' + \bar\partial$ be the Chern connection of $(T^{1,\,0}X,\,\omega)$. 

From $\bar\partial\alpha_j = \alpha_j\wedge\overline\alpha_j/2$ for $j=1,\dots , s$, we deduce that $$\bar\partial e_j = \frac{\overline\alpha_j}{2}\otimes e_j, \hspace{5ex} j=1,\dots , s.$$ On the other hand, from $\langle e_j,\,e_k\rangle_\omega = \delta_{jk}$ (the Kronecker delta) we deduce the first equality below: $$0 = D' \langle e_j,\,e_k\rangle_\omega = \{D'e_j,\,e_k\} + \{e_j,\,\bar\partial e_k\} = \{D'e_j,\,e_k\} + \delta_{jk}\,\frac{\alpha_k}{2}, \hspace{5ex} j,\,k=1,\dots , s,$$ where the sesquilinear bracket \begin{eqnarray*}\{\cdot,\,\cdot\}:C^\infty_p(X,\,T^{1,\,0}X)\times C^\infty_q(X,\,T^{1,\,0}X) & \longrightarrow & C^\infty_{p+q}(X,\,\C), \\
  \bigg\{\sum\limits_\lambda\sigma_\lambda\otimes g_\lambda,\,\sum\limits_\mu\tau_\mu\otimes g_\mu\bigg\} & := & \sum\limits_{\lambda,\,\mu}\sigma_\lambda\wedge\overline\tau_\mu\,\langle g_\lambda,\,g_\mu\rangle_\omega\end{eqnarray*} combines the wedge product of scalar-valued forms (in this case, the locally defined forms $\sigma_\lambda$ and $\tau_\mu$ that represent the given $T^{1,\,0}X$-valued forms in a local trivialisation of $T^{1,\,0}X$ defined by a local frame $(g_\lambda)_\lambda$) with the inner product defined by the Hermitian metric $\omega$ on the fibres of $T^{1,\,0}X$ (see [Dem97, V.7.]).

Hence, we get: $$D'e_j = -\frac{\alpha_j}{2}\otimes e_j, \hspace{5ex} j=1,\dots , s.$$

Letting $\Theta:=\Theta_\omega(T^{1,\,0}X) = D^2 = D'\bar\partial + \bar\partial D'$ be the Chern curvature form of $(T^{1,\,0}X, \omega)$, we get: \begin{eqnarray*}\Theta e_j & = & D'\bigg(\frac{\overline\alpha_j}{2}\otimes e_j\bigg) + \bar\partial\bigg(-\frac{\alpha_j}{2}\otimes e_j\bigg) = \frac{1}{2}\,\partial\overline\alpha_j\otimes e_j + \frac{1}{2}\,\overline\alpha_j\otimes D'e_j - \frac{1}{2}\,(\bar\partial\alpha_j)\otimes e_j - \frac{1}{2}\,\alpha_j\otimes\bar\partial e_j \\
 & = & -\frac{1}{4}\,(\alpha_j\wedge\overline\alpha_j)\otimes e_j -\frac{1}{4}\,(\overline\alpha_j\wedge\alpha_j)\otimes e_j -\frac{1}{4}\,(\alpha_j\wedge\overline\alpha_j)\otimes e_j -\frac{1}{4}\,(\alpha_j\wedge\overline\alpha_j)\otimes e_j \\
& = & -\frac{1}{2}\,(\alpha_j\wedge\overline\alpha_j)\otimes e_j, \hspace{5ex} j=1,\dots , s.\end{eqnarray*}

It follows that the holomorphic sectional curvature in the direction of the tangent vector $e_j$ is $$\bigg\langle i\Theta(e_j,\,\bar{e}_j)e_j,\,e_j\bigg\rangle_\omega = -\frac{1}{2}\,\bigg\langle(i\alpha_j\wedge\overline\alpha_j)(e_j,\,\bar{e}_j)e_j,\,e_j\bigg\rangle_\omega = -\frac{1}{2}\,\langle e_j,\,e_j\rangle_\omega = -\frac{1}{2}$$ for every $j=1,\dots , s$. 

We have thus proved the following

\begin{Prop}\label{Prop:O-T_hol-sectional_E}
For every  solvmanifold $X=G/\Lambda$ as in Proposition \ref{Prop:O-T_partially-s-K_hyp} (e.g. O-T manifold), the {\bf holomorphic sectional curvatures} of $X$ in the directions of the horizontal vector bundle $E$ are constant, equal to $-\frac{1}{2}$. 

\end{Prop}

\subsubsection{Curvature of the canonical bundle}\label{curvature-K_X_O-T} Starting from the definitions (\ref{eqn:O-T_alpha_j-beta_k}) of the $C^\infty$ $(1,\,0)$-forms $\alpha_j$ and $\beta_k$, we consider the real-valued functions $$X_j:=e^{x_j}, \hspace{5ex} j=1,\dots , s$$ and the complex-valued functions $$Z_j:=X_j + iy_j, \hspace{5ex} j=1,\dots , s.$$ 

\begin{Claim}\label{Claim:O-T_local-hol-frame} $\{dZ_1, \dots , dZ_s,\,dz_1,\dots , dz_t\}$ is a local {\it holomorphic frame} of the holomorphic cotangent bundle $\Lambda^{1,\,0}T^\star X$.

\end{Claim}

\noindent {\it Proof.} From (\ref{eqn:O-T_alpha_j-beta_k}) we get: \begin{eqnarray*}\alpha_j = d\log X_j + \frac{i}{X_j}\,dy_j = \frac{1}{X_j}\,(dX_j + idy_j)= \frac{1}{X_j}\,dZ_j, \hspace{5ex} j=1,\dots , s,\end{eqnarray*} hence $dZ_j = X_j\,\alpha_j$ is indeed a $(1,\,0)$-form for every $j=1,\dots , s$. Since $z_1,\dots , z_t$ are the complex variables of $\C^t$, the $1$-forms $dz_1,\dots dz_t$ are of type $(1,\,0)$ and holomorphic. To see that the $(1,\,0)$-forms $dZ_1,\dots , dZ_s$ are holomorphic, we compute: \begin{eqnarray*}\bar\partial(dZ_j) = \bar\partial(e^{x_j}\,\alpha_j) = e^{x_j}\,\bar\partial x_j\wedge\alpha_j + e^{x_j}\,\bar\partial\alpha_j = \frac{e^{x_j}}{2}\,\bar\alpha_j\wedge\alpha_j + e^{x_j}\,\frac{\alpha_j\wedge\bar\alpha_j}{2} = 0,\end{eqnarray*} where we used the equality $\bar\partial x_j = (1/2)\,\bar\alpha_j$ that follows from $$dx_j = \mbox{Re}(\alpha_j) = \frac{\alpha_j + \bar\alpha_j}{2},$$ from $\bar\partial x_j$ being the $(0,\,1)$-part of $dx_j$ and from $\alpha_j$ being of type $(1,\,0)$. \hfill $\Box$

\vspace{2ex}

The inner products of the pairs of elements of the above frame are: \begin{eqnarray*}\langle dZ_j,\,dZ_k\rangle_\omega & = & X_jX_k\,\langle\alpha_j,\,\alpha_k\rangle_\omega = X_jX_k\,\delta_{jk}, \hspace{5ex} j,k=1,\dots , s,\\
 \langle dz_l,\,dz_r\rangle_\omega & = & e^{\psi_l + \overline\psi_r}\,\langle\beta_l,\,\beta_r\rangle_\omega = e^{\psi_l + \overline\psi_r}\,\delta_{lr}, \hspace{5ex} l,r=1,\dots , t, \\
 \langle dZ_j,\,dz_r\rangle_\omega & = & X_je^{\overline\psi_r}\,\langle\alpha_j,\,\beta_r\rangle_\omega = 0, \hspace{5ex} j=1,\dots , s \hspace{1ex}\mbox{and}\hspace{1ex}  r=1,\dots , t.\end{eqnarray*}

Thus, for the non-vanishing local holomorphic frame of the canonical bundle $K_X$ of the solvmanifold $X$ given by the $(n,\,0)$-form $e=dZ_1\wedge\dots\wedge dZ_s\wedge dz_1\wedge\dots\wedge dz_t$, the squared pointwise $\omega$-norm is $|e|^2_\omega = X_1^2\dots X_s^2\,e^{2\mbox{Re}\,\psi_1 +\dots+ 2\mbox{Re}\,\psi_t}.$ Hence, the weight function $\varphi$ defined by $e^{-\varphi} = |e|_\omega$ is $\varphi = -\sum_{j=1}^sx_j - \sum_{l=1}^t\mbox{Re}\,\psi_l$. Now, by the unimodularity of the Lie group $G$, the determinant of the matrix $$\mbox{diag}\bigg(e^{x_1}, \dots , e^{x_s}, e^{\psi_1}, \dots , e^{\psi_t},\,e^{\overline\psi_1},\dots , e^{\overline\psi_t} \bigg)$$ equals $1$, which translates to $\sum_{l=1}^t\mbox{Re}\,\psi_l = -(1/2)\,\sum_{j=1}^sx_j$. This leads to $\varphi = -(1/2)\,\sum_{j=1}^sx_j$. 

The curvature form of the canonical bundle $K_X$ of the solvmanifold $X$ is then given by $$i\Theta_\omega(K_X) = i\partial\bar\partial\varphi = -\frac{1}{2}\,\sum_{j=1}^si\partial\bar\partial x_j.$$ On the other hand, as observed above, $\bar\partial x_j = (1/2)\,\bar\alpha_j$, so we get $$\partial\bar\partial x_j = \frac{1}{2}\partial\overline\alpha_j = \frac{1}{2}\,\overline{d\alpha_j} = -\frac{1}{4}\,\alpha_j\wedge\overline\alpha_j, \hspace{5ex} j=1,\dots , s,$$ where we used the equalities $\bar\partial\alpha_j = d\alpha_j = (1/2)\,\alpha_j\wedge\overline\alpha_j$. 

The conclusion of this computation is summed up in the following

\begin{Prop}\label{Prop:O-T_curvature_K_X}For every  solvmanifold $X=G/\Lambda$ as in Proposition \ref{Prop:O-T_partially-s-K_hyp}(e.g. O-T manifold), the {\bf curvature form} of the canonical bundle $K_X$ with respect to the fibre metric induced by the Hermitian $\omega$ defined in (\ref{eqn:O-T_metric}) is \begin{eqnarray*}i\Theta_\omega(K_X) = \frac{1}{8}\,\sum\limits_{j=1}^si\alpha_j\wedge\overline\alpha_j.\end{eqnarray*} 

In particular, $i\Theta_\omega(K_X)\geq 0$ on $X$ and $i\Theta_\omega(K_X)(x)(\xi,\,\bar\xi)> 0$ for every $x\in X$ and every $\xi\in E_x$.

\end{Prop}

\subsection{Partial hyperbolicity of Miebach-Oeljeklaus manifolds}\label{subsection:examples:M-O} 
Consider the matrix group
\[H=\left\{\left(
\begin{array}{ccc}
1& y &v \\
0&   x&z\\
0&0&1  
\end{array}
\right): x,y,z,v\in \R\right\}.\]
Then, the Lie algebra $\frak h$ of $H$ is spanned by
\[A=\left(
\begin{array}{ccc}
0& 0 &0 \\
0&   1&0\\
0&0&0 
\end{array}
\right), B=\left(
\begin{array}{ccc}
0& 1 &0 \\
0&   0&0\\
0&0&0 
\end{array}
\right), C=\left(
\begin{array}{ccc}
0& 0 &0 \\
0&   0&1\\
0&0&0 
\end{array}
\right), T=\left(
\begin{array}{ccc}
0& 0 &1 \\
0&   0&0\\
0&0&0 
\end{array}
\right).
\]
This leads to the structure equations:
\[[A,B]=-B, \hspace{2ex} [A,C]=C, \hspace{2ex} [B,C]=T, \hspace{2ex} [A,T]=[B,T]=[C,T]=0.
\]

Now, take the dual basis $\{a, b,c,t\}$ of $\frak h^{\ast}$. Then $da=0, db=a\wedge b, dc=-a\wedge c, d t=-b\wedge c$. Define $\alpha=a+ib$ and $\beta=c+it$. We get: $$d\alpha=ia\wedge b =-\frac{\alpha\wedge\bar\alpha}{2} \hspace{3ex} \mbox{and}  \hspace{3ex} d\beta=- (a+ib)\wedge c=-\frac{\alpha\wedge (\beta+\bar\beta)}{2}.$$ Regarding $\{\alpha, \beta\}$ as a global $\C^{\infty}$-frame of $(1,0)$-forms, we have a left-invariant complex structure on $H$.

Now, consider $G=H^d$. Miebach and Oeljeklaus construct in [MO22] a lattice $\Gamma$ in $G$ corresponding to a totally real number field $K$ of degree $d$ satisfying certain conditions. From this, we get the compact complex manifold $X:=G/\Gamma$ with a global $C^{\infty}$-frame $\{\alpha_{1},\dots, \alpha_{d},\,\beta_{1},\dots, \beta_{d}\}$ of $(1,0)$-forms satisfying the following equations: $$d\alpha_{i}=-\frac{\alpha_{i}\wedge\bar\alpha_{i}}{2} \hspace{3ex} \mbox{and} \hspace{3ex} d\beta_i=-\frac{\alpha_{i}\wedge (\beta_{i}+\bar\beta_{i})}{2}$$ for all $i\in\{1,\dots , d\}$. Thus, $\bar\partial \alpha_{i}=-\frac{\alpha_{i}\wedge\bar\alpha_{i}}{2}$ and $\bar\partial\beta_{i}=-\frac{\alpha_{i}\wedge \bar\beta_{i}}{2}$.
Let  $\{E_{1},\dots, E_{d},\, F_{1},\dots, F_{d}\}$ be the frame of $T^{1,0}X$ dual to the frame $\{\alpha_{1},\dots, \alpha_{d},\,\beta_{1},\dots, \beta_{d}\}$ of $\Lambda^{1,\,0}T^\star X$. We get:  $$\bar\partial E_{i}=-\frac{\bar\alpha_{i}\wedge E_{i}+\bar\beta_{i}\wedge F_{i}}{2} \hspace{3ex} \mbox{and} \hspace{3ex} \bar\partial F_{i}=0$$ for all $i\in\{1,\dots , d\}$.

\vspace{2ex}

Let $E$ and $F$ be the $C^\infty$-trivial vector subbundles of $T^{1,\,0}X$ generated by $E_1,\dots , E_d$ and, respectively, by $F_1,\dots , F_d$. Note that while $F$ is a holomorphic subbundle of $T^{1,\,0}X$, $E$ is only a $C^\infty$ subbundle. 

A natural Hermitian metric on $X$ is defined by the positive definite $C^\infty$ $(1,\,1)$-form \begin{eqnarray}\label{eqn:O-T_metric}\omega = i\alpha_1\wedge\overline\alpha_1 + \dots + i\alpha_d\wedge\overline\alpha_d + i\beta_1\wedge\overline\beta_1 + \dots + i\beta_d\wedge\overline\beta_d>0.\end{eqnarray} 
Thus, $\omega = \omega_E + \omega_F$, where $\omega_E:=i\alpha_1\wedge\overline\alpha_1 + \dots + i\alpha_d\wedge\overline\alpha_d\geq 0$ and $\omega_F:=i\beta_1\wedge\overline\beta_1 + \dots + i\beta_d\wedge\overline\beta_td\geq 0$ are $C^\infty$ positive semi-definite $(1,\,1)$-forms on $X$. By construction, $\omega_E(\xi,\,\bar\xi)>0$ for every $(1,\,0)$-tangent vector $\xi$ of $E$ and $\omega_E(\eta,\,\bar\nu) = 0$ for all $(1,\,0)$-tangent vectors $\eta,\,\nu$ of $F$. The analogous property is satisfied by $\omega_F$ with $E$ and $F$ permuted.

Let us now consider the following $C^\infty$ $(1,\,1)$-form on $X$: \begin{eqnarray*}\Omega = \omega_E= d\bigg(-2i\alpha_1-\cdots -2i\alpha_d\bigg).\end{eqnarray*} Since $\Omega$ is $d$-exact on the {\it compact} complex manifold $X$, its lift to the universal cover $G$ is $d(bounded)$ with respect to the lift of the metric $\omega$. In particular, the triple $(E,\,\Omega,\,\omega= \omega_E + \omega_F)$ is a partially $1$-K\"ahler hyperbolic structure on $X$. 

We have thus proved the following

\begin{Prop}\label{Prop:M-O} The solvmanifold $X=G/\Lambda$ obtained as the quotient of $G=H^d$ by any co-compact lattice $\Lambda$ (e.g. the Miebach-Oeljeklaus manifold) is {\bf partially $1$-K\"ahler hyperbolic}. 
\end{Prop}

\subsection{Partial hyperbolicity of a class of complex parallelisable solvmanifolds}\label{subsection:class_solv-Cpar}

Let $n\in\N^\star$ and let $G=\C^n\ltimes\C^{n+1}$ be the complex Lie group defined as a semi-direct product via the map \begin{eqnarray*}\C^n\ni(z_1,\dots , z_n)\mapsto\mbox{diag}\bigg(e^{z_1},\dots , e^{z_n},\,e^{-z_1-\dots - z_n}\bigg)\in\mbox{Aut}(\C^{n+1})=GL_{n+1}(\C),\end{eqnarray*} where $\mbox{diag}$ stands for the diagonal matrix with the indicated diagonal entries. This Lie group is solvable.

We can construct a lattice $\Gamma$ in $G$ associated with a totally real algebraic number field K of degree $n + 1$ in the same way as in Example 4 of [Kas17]. To wit, consider the semi-direct product  $H=\R^{n}\ltimes_{\psi} \R^{n+1}$ via the map 
\[\R^n\ni (x_{1},\dots, x_{n})\mapsto\mbox{diag}(e^{x_{1}},\dots e^{x_{n}}, e^{-x_{1}-\dots-x_{n}})\in GL_{n+1}(\R).
\]
For a totally real  algebraic number field $K$ of degree $n+1$, Dirichlet's unit theorem yields a subgroup $\Gamma^{\prime}\subset {\mathcal O}_{K}^{\ast}$ such that $\Gamma^{\prime}$ can be regarded as a lattice in $\R^{n}$. Then, the semi-direct product $ \Gamma^{\prime}_{1}\ltimes {\mathcal O}_{K}$ can be regarded as a lattice in  $H=\R^{n}\ltimes_{\psi} \R^{n+1}$ (see [TY09]). We obtain the lattice $$\Gamma=(\Gamma^{\prime}_{1}\oplus (2\pi\Z)^n)\ltimes ({\mathcal O}_{K}+\sqrt{-1}{\mathcal O}_{K})$$ in $G$. We will investigate the compact complex manifold $X=G/\Gamma$, a solvmanifold of complex dimension $2n+1$.

Denoting by $z_1,\dots , z_n$ the coordinates of $\C^n$ and by $w_1,\dots , w_{n+1}$ the coordinates of $\C^{n+1}$, the $C^\infty$ $(1,\,1)$-form \begin{eqnarray}\label{eqn:G-metric_def}\widetilde\omega = \sum\limits_{j=1}^n idz_j\wedge d\bar{z}_j + \sum\limits_{j=1}^n e^{-z_j-\bar{z}_j}\,idw_j\wedge d\overline{w}_j + e^{\sum\limits_{j=1}^n(z_j+\bar{z}_j)}\,idw_{n+1}\wedge d\overline{w}_{n+1}\end{eqnarray} is positive definite, hence it defines a Hermitian metric, on $G$. Moreover, it passes to the quotient and induces the Hermitian metric $$\omega=\sum\limits_{j=1}^n i\alpha_j\wedge\overline\alpha_j + \sum\limits_{j=1}^n i\beta_j\wedge\overline\beta_j + i\beta_{n+1}\wedge\overline\beta_{n+1}$$ on $X$, where $\alpha_1,\dots , \alpha_n$, $\beta_1,\dots , \beta_n$ and $\beta_{n+1}$ are the holomorphic $(1,\,0)$-forms induced on $X$ by the holomorphic $(1,\,0)$-forms $dz_1,\dots , dz_n$, $e^{-z_1}\,dw_1,\dots , e^{-z_n}\,dw_n$ and, respectively, $e^{\sum\limits_{j=1}^nz_j}\,dw_{n+1}$ of $G$. 

Since $G$ is a {\it complex} Lie group, the manifold $X$ is {\it complex parallelisable} in the sense that its holomorphic tangent bundle $T^{1,\,0}X$ is holomorphically trivial. A global holomorphic frame for the holomorphic cotangent bundle $\Lambda^{1,\,0}T^\star X$ is provided by $\{\alpha_1,\dots , \alpha_n,\,\beta_1,\dots ,\beta_n,\beta_{n+1}\}$. Let $\{e_1,\dots , e_n,\,f_1,\dots , f_n,f_{n+1}\}$ be the dual global holomorphic frame of $T^{1,\,0}X$. 

It is easy to see that the $(2n,\,2n)$-form $\widetilde\omega_{2n}:=\widetilde\omega^{2n}/(2n)!$ is $d$-exact on $G$ (the universal cover of $X$), but it is not $d(bounded)$ on $G$. However, we will now show that if we remove one of the terms $\exp(-z_j-\bar{z}_j)\,idw_j\wedge d\overline{w}_j$ (for example, the one corresponding to $j=1$) from the sum defining $\widetilde\omega$, we get a semi-positive $(1,\,1)$-form whose $(2n)$-th power is $d(bounded)$ on $G$. This will lead to a {\it partially $(2n)$-K\"ahler hyperbolic} structure on $X$.

Let $$\widetilde\omega_1 = \sum\limits_{j=1}^n idz_j\wedge d\bar{z}_j + \sum\limits_{j=2}^n e^{-z_j-\bar{z}_j}\,idw_j\wedge d\overline{w}_j + e^{\sum\limits_{j=1}^n(z_j+\bar{z}_j)}\,idw_{n+1}\wedge d\overline{w}_{n+1}$$ be the positive semi-definite $C^\infty$ $(1,\,1)$-form on $G$ equal to $\widetilde\omega - e^{-z_1-\bar{z}_1}\,idw_1\wedge d\overline{w}_1$. We have: \begin{eqnarray*}(\widetilde\omega_1)_{2n} & = & \frac{\omega_1^{2n}}{(2n)!} = e^{z_1+\bar{z}_1}\,idz_1\wedge d\bar{z}_1\wedge\prod\limits_{l=2}^nidz_l\wedge d\bar{z}_l\wedge\bigg(\prod\limits_{r=2}^nidw_r\wedge d\overline{w}_r\bigg)\wedge(idw_{n+1}\wedge d\overline{w}_{n+1})\\
 & = & d\bigg(ie^{z_1+\bar{z}_1}\,d\bar{z}_1\wedge\prod\limits_{l=2}^nidz_l\wedge d\bar{z}_l\wedge\prod\limits_{r=2}^nidw_r\wedge d\overline{w}_r\wedge idw_{n+1}\wedge d\overline{w}_{n+1}\bigg) \\
 & = & d\bigg(i\,d\bar{z}_1\wedge\prod\limits_{l=2}^nidz_l\wedge d\bar{z}_l\wedge\prod\limits_{r=2}^n\bigg(e^{-z_r-\bar{z}_r}\,idw_r\wedge d\overline{w}_r\bigg)\wedge\bigg(e^{\sum\limits_{j=1}^n(z_j+\bar{z}_j)}\,idw_{n+1}\wedge d\overline{w}_{n+1}\bigg)\bigg)\\
 & = & d\widetilde\Gamma,\end{eqnarray*} with the definition of the $(2n-1)$-form $\widetilde\Gamma$ on $G$ made obvious by the notation.

A comparison of the formula for $\widetilde\Gamma$ with the definition (\ref{eqn:G-metric_def}) of $\widetilde\omega$ shows that $\widetilde\Gamma$ is $\widetilde\omega$-bounded on $G$. We conclude that $(\widetilde\omega_1)_{2n}$ is $d(bounded)$ on $G$ w.r.t. the metric $\widetilde\omega$, or equivalently, that the $(2n,\,2n)$-form $(\omega_1)_{2n}$ is $\widetilde{d}(bounded)$ on $X$ w.r.t. the metric $\omega$, where $\omega_1$ is the positive semi-definite $C^\infty$ $(1,\,1)$-form $$\omega_1=\sum\limits_{j=1}^n i\alpha_j\wedge\overline\alpha_j + \sum\limits_{j=2}^n i\beta_j\wedge\overline\beta_j + i\beta_{n+1}\wedge\overline\beta_{n+1}$$ on $X$. This shows that $(E,\,\Omega,\,\omega= \omega_E + \omega_{nE})$ is a {\it partially (2n)-K\"ahler hyperbolic structure} on $X$, where $E$ is the (globally trivial) holomorphic subbundle of $T^{1,\,0}X$ generated by $e_1,\dots , e_n,\,f_2,\dots , f_n, f_{n+1}$ and $\Omega:=(\omega_1)_{2n}$. We put, of course, $\omega_E:= \omega_1$ and $\omega_{nE}:=i\beta_1\wedge\overline\beta_1$.   

We have thus proved the following

 \begin{Prop}\label{Prop:C-par-solvmanifolds} The $(2n+1)$-dimensional complex parallelisable solvmanifolds $X=G/\Gamma$ described above are {\bf partially (2n)-K\"ahler hyperbolic}.

\end{Prop}

Since the metric $\omega$ is flat (as can be easily seen) on the globally holomorphically trivial vector bundle $T^{1,\,0}X$, the induced holomorphic sectional curvatures of $X$ and the curvature of $K_X$ vanish.

\subsection{Vaisman manifolds}\label{subsection:Vaisman} Let $(X,J)$ be a compact complex manifold endowed with a Hermitian metric $g$ and let $\omega=g(\cdot,\,J\cdot)$ be the fundamental form of $g$. The following definition is standard: the metric $g$ is said to be {\it locally conformal K\"ahler} (lcK) if there exists a closed $1$-form $\theta$ (called the Lee form) such that $d\omega=\theta\wedge \omega$.

It is known (see Theorem 2.1 of [DO98]) that if $\theta$ is not exact, the manifold $(X,J)$ does not admit any K\"ahler structure.

Let $\nabla$ be the Levi-Civita connection of $g$. An lcK metric $g$ is said to be a {\it Vaisman metric} ([Vai79]) if  $\nabla \theta=0$.

Now, suppose that $g$ is Vaisman. We summarise below the argument in Sections 2 and 3  of [Tsu94]. Let $A$ and $B$ be the vector fields that are dual to the $1$-forms $\theta$ and respectively $-\theta \circ J$ with respect to the metric $g$. Then $$A=JB, \hspace{3ex} L_{A}J=0, \hspace{3ex} L_{B}J=0, \hspace{3ex} L_{A}g=0, \hspace{3ex} L_{B}g=0 \hspace{3ex} \mbox{and} \hspace{3ex} [A,B]=0.$$ 

The holomorphic vector field $B-iA$ generates a holomorphic foliation $\mathcal F$. We have $$\omega=d(\theta\circ J)-\theta\wedge(\theta \circ J)$$ and $d(\theta\circ J)$ is a transverse K\"ahler structure on $\mathcal F$. This means that $$d(\theta\circ J)=0 \hspace{1ex} \mbox{on} \hspace{1ex} \mathcal F \hspace{3ex} \mbox{and} \hspace{3ex} d(\theta\circ J)>0 \hspace{1ex} \mbox{on} \hspace{1ex} T^{1,0}/\mathcal F.$$ Consider the $C^{\infty}$ vector subbundle $E=ker(\theta + i\,\theta \circ J)$ of the holomorphic  tangent vector bundle $ T^{1,0}X$ of $X$. 

Since the $(1,1)$-form $\Omega=d(\theta\circ J)$ is $d$-exact on the {\it compact} complex manifold $X$, its lift to the universal cover is $d(bounded)$ with respect to the lift of the metric $\omega$. Hence, $(E,\,\Omega,\,\omega=d(\theta\circ J)-\theta\wedge(\theta \circ J))$ is a partially $1$-K\"ahler hyperbolic structure on $X$. We have thus proved the following

\begin{Prop}\label{Prop:Vaisman} Every compact complex manifold admitting a Vaisman metric is {\bf partially $1$-K\"ahler hyperbolic}.  

\end{Prop}
We remark that a compact Vaisman manifold is never Kobayashi   hyperbolic.

A primary Hopf manifold is a compact complex manifold obtained as the quotient of $\C^n\setminus\{0\}$ by a subgroup generated by a transformation of $\C^n\setminus\{0\}$ of the shape $$(z_1, \dots , z_n) \mapsto (\lambda_1z_1,\dots , \lambda_n z_n),$$ where $\lambda_1,\dots , \lambda_n$ are complex numbers such that $0<\vert \lambda_{n}\vert \le  \dots\le \vert \lambda_{1}\vert<1$.

It is known (see [KO05]) that every primary Hopf manifold admits a Vaisman metric. Thus, every primary Hopf manifold is partially $1$-K\"ahler hyperbolic.

\section{Ahlfors currents}\label{section:Ahlfors} In this section, we give a sufficient condition, reminiscent of de Th\'elin's criteria of [dT10], for the existence of an Ahlfors current on a compact complex Hermitian manifold. However, our condition seems simpler, is cast in the language of this paper and is similar to the {\it subexponential growth condition} of Definition \ref{Def:subexp}. Thus, it seems better suited to our situation. We then go on to discuss several examples and the link between these currents and partial hyperbolicity.

\subsection{Existence of Ahlfors currents}\label{subsection:Ahlfors_existence}

Let $f:\C^p\longrightarrow (X,\,\omega)$ be a holomorphic map that is non-degenerate at some point $x_0\in\C^p$, where $1\leq p\leq n-1$ and $X$ is an $n$-dimensional compact complex manifold equipped with a Hermitian metric $\omega$. We will use the notation of $\S.$\ref{subsection:partially-p-hyperbolic} (with $q=p$). In particular, for every $r>0$, $B_r$ and $S_r$ stand for the open ball, respectively the sphere, of radius $r$ centred at the origin of $\C^p$, while $[B_r]$ and $[S_r]$ denote the currents of integration thereon. One can consider the direct image $f_\star[B_r]$ of the current $[B_r]$ under $f$ (cf. e.g. [Dem97, I-2.C.1.]). It is a current in $X$ of the same bidimension $(p,\,p)$ as $[B_r]$, so $f_\star[B_r]\geq 0$ is a strongly positive current of bidegree $(n-p,\,n-p)$ in $X$. (See e.g. [Dem97, III-1.C.] for the notions of strongly positive and (weakly) positive currents.) 

We can normalise $f_\star[B_r]$ to get the bidegree-$(n-p,\,n-p)$-current \begin{eqnarray}\label{eqn:T_r_def}T_r:=\frac{1}{\mbox{Vol}_{\omega,\,f}(B_r)}\,f_\star[B_r]\geq 0,  \hspace{6ex} r>0,\end{eqnarray} in $X$. It has unit mass with respect to $\omega$: \begin{eqnarray}\label{eqn:unit-mass_T-r}\int\limits_XT_r\wedge\omega_p = \frac{1}{\mbox{Vol}_{\omega,\,f}(B_r)}\,\int\limits_{\C^p}[B_r]\wedge f^\star\omega_p = \frac{1}{\mbox{Vol}_{\omega,\,f}(B_r)}\,\int\limits_{B_r}f^\star\omega_p = 1,  \hspace{6ex} r>0.\end{eqnarray} Thus, the family $(T_r)_{r>0}$ of strongly positive currents, being uniformly bounded in mass, has a weakly convergent subsequence $(T_{r_\nu})_{\nu\in\N}$ with $r_\nu\to +\infty$. The limiting current $T\geq 0$ is strongly positive of bidegree $(n-p,\,n-p)$ in $X$. However, $T$ need not be $d$-closed since \begin{eqnarray}\label{eqn:Stokes_currents}d[B_r] = -[S_r]\neq 0,  \hspace{6ex} r>0,\end{eqnarray} as follows at once from the Stokes theorem. 

A current $T$ obtained as the limit in the weak topology of currents of a sequence of currents $T_{r_\nu}\geq 0$ constructed as above from a holomorphic map $f:\C^p\longrightarrow (X,\,\omega)$, with $r_\nu\to +\infty$, is said to be an {\bf Ahlfors current} if $dT=0$. It is well known that Ahlfors currents need not exist on an arbitrary compact complex manifold $X$.

 \vspace{2ex}

Before giving our sufficient condition for the existence of an Ahlfors current $T$, we give a general estimate on the norms of the currents $\partial T_r$ that holds without any special assumption on the map $f$.

 We start by following de Th\'elin's strategy of [dT10] from which we will deviate at some point that will be specified. Since $T$ is the limit of strongly positive (hence real) currents, it is itself real. Thus, $\bar\partial T$ is the conjugate of $\partial T$, so proving that $dT=0$ is equivalent to proving that $\partial T = 0$. To this end, we will show that a certain norm of $\partial T$ vanishes. 

As a current of bidegree $(n-p+1,\,n-p)$ on the $n$-dimensional compact complex manifold $X$, $\partial T$ acts on $C^\infty$ forms of bidegree $(p-1,\,p)$. We consider the closed unit ball of these forms with respect to the $C^0$-norm induced by the metric $\omega$: $${\cal F}_\omega(p-1,\,p):=\{\psi\in C^\infty_{p-1,\,p}(X,\,\C)\,\mid\,||\psi||_{C^0_\omega}:=\max\limits_{x\in X}|\psi(x)|_\omega\leq 1\}$$ and the induced norm on the space ${\cal D}^{'n-p+1,\,n-p}(X)$ of bidegree-$(n-p+1,\,n-p)$-currents $S$ on $X$: $$||S||:=\sup\limits_{\psi\in{\cal F}_\omega(p-1,\,p)}|\langle S,\,\psi\rangle|.$$

Our general estimate, different from the one in [dT10, Theorem 0.1], is spelt out in the following 

\begin{The}\label{The:del-T_r-norm_estimate} Let $X$ be an $n$-dimensional compact complex manifold equipped with a Hermitian metric $\omega$. Let $p\in\{1,\dots , n-1\}$ and let $f:\C^p\longrightarrow X$ be a holomorphic map that is non-degenerate at some point $x_0\in\C^p$.

Then, the norm of $\partial T_r$, where $T_r$ is the current defined by $f$ through formula (\ref{eqn:T_r_def}), satisfies the estimate: \begin{eqnarray}\label{eqn:T_r_norm-estimate}||\partial T_r||\leq\frac{1}{\sqrt{2}}\,\frac{A_{\omega,\,f}(S_r)}{\mbox{Vol}_{\omega,\,f}(B_r)},  \hspace{6ex} r>0,\end{eqnarray} where $A_{\omega,\,f}(S_r)$ and $\mbox{Vol}_{\omega,\,f}(B_r)$ are the $(\omega,\,f)$-area of the Euclidean sphere $S_r\subset\C^p$, respectively the $(\omega,\,f)$-volume of the Euclidean ball $B_r\subset\C^p$, defined in (\ref{eqn:omega-f_area_def}) and (\ref{eqn:omega-f_vol_def}). 

\end{The}

\noindent {\it Proof.} $\bullet$ In a departure from the strategy of [dT10], we handle the $(p-1,\,p)$-forms on $\C^p$ starting from the standard fact that the pointwise multiplication by the $(p-1)$-st power of any metric (in our case, of the degenerate metric $f^\star\omega$, that is a genuine metric on $\C^p\setminus\Sigma_f$) is an isomorphism on its image when acting on $1$-forms in $\C^p$. In particular, the linear map: \begin{eqnarray*}f^\star\omega_{p-1}\wedge\cdot:\Lambda^{0,\,1}T^\star\C^p\longrightarrow\Lambda^{p-1,\,p}T^\star\C^p\end{eqnarray*} is bijective at every point of $\C^p\setminus\Sigma_f$. Hence, for every $(p-1,\,p)$-form $\psi$ on $X$, there exists a unique $(0,\,1)$-form $\alpha_\psi$ on $\C^p$ such that \begin{eqnarray}\label{eqn:f-star-psi_splitting}f^\star\psi = \alpha_\psi\wedge f^\star\omega_{p-1}.\end{eqnarray}

This gives the latter equality below, while the former equality follows from the standard formula (\ref{eqn:prim-form-star-formula-gen}) expressing the image under the Hodge star operator of any primitive form (in this case, $\alpha_\psi$, which is a $1$-form, hence primitive): \begin{eqnarray*}\label{eqn:f-star-psi_Hodge-star}\star_{f^\star\omega}\alpha_\psi = i\,\alpha_\psi\wedge f^\star\omega_{p-1} = i\, f^\star\psi.\end{eqnarray*} 

In particular, since the Hodge star operator is an isometry w.r.t. the pointwise norm, we infer \begin{eqnarray}\label{eqn:f-star-psi_Hodge-star_norms}|\alpha_\psi|_{f^\star\omega} = |f^\star\psi|_{f^\star\omega} = |\psi|_\omega\end{eqnarray} at every point in $\C^p\setminus\Sigma_f$.

\vspace{1ex}

$\bullet$ We now set about proving estimate (\ref{eqn:T_r_norm-estimate}). For $\psi\in C^\infty_{p-1,\,p}(X,\,\C)$ and $r>0$, we have: \begin{eqnarray}\label{eqn:proof_main-estimate_1}\nonumber\langle\partial f_\star[B_r],\,\psi\rangle & = & \langle\partial[B_r],\,f^\star\psi\rangle \stackrel{(a)}{=}\langle d[B_r],\,f^\star\psi\rangle \stackrel{(b)}{=} -\langle[S_r],\,\alpha_\psi\wedge f^\star\omega_{p-1}\rangle \\
  & = & -\int\limits_{S_r}\alpha_\psi\wedge f^\star\omega_{p-1}:\stackrel{(c)}{=}h(r) \stackrel{(d)}{=} H'(r),\end{eqnarray} where (a) follows from $\langle\bar\partial[B_r],\,f^\star\psi\rangle = 0$ (which holds trivially, for bidegree reasons, since $\bar\partial[B_r]$ is a current of bidegree $(0,\,1)$ and $f^\star\psi$ is a form of bidegree $(p-1,\,p)$ in $\C^p$), (b) follows from (\ref{eqn:Stokes_currents}) and (\ref{eqn:f-star-psi_splitting}), (c) is the definition of a function $h:(0,\,+\infty)\to\R$ and (d) is the immediate consequence of the definition $H(r):=\int\limits_0^rh(t)\,dt$ of a function $H:(0,\,+\infty)\to\R$.

As in [dT10], we momentarily fix arbitrary reals $0<r<r'$, but we will handle the integrals differently. We have: \begin{eqnarray*}\frac{H(r') - H(r)}{r'-r} = -\frac{1}{r'-r}\,\int\limits_r^{r'}\bigg(\int\limits_{S_t}\alpha_\psi\wedge f^\star\omega_{p-1}\bigg)\,dt = -\frac{1}{r'-r}\,\int\limits_{B_{r'}\setminus B_r}\alpha_\psi\wedge d\rho\wedge f^\star\omega_{p-1},\end{eqnarray*} where the last equality follows from the Fubini theorem and $\rho:\C^p\to[0,\,+\infty)$ is the function defined by $\rho(z)=|z|$ (hence $\rho = t$ on $S_t$ and $d\rho = dt$). Since $d\rho = \partial\rho + \bar\partial\rho$ and $\alpha_\psi\wedge\bar\partial\rho\wedge f^\star\omega_{p-1} = 0$ for bidegree reasons (as it is a $(p-1,\,p+1)$-form in $\C^p$), we get the first equality below: \begin{eqnarray}\label{eqn:H-estimate_1}\frac{H(r') - H(r)}{r'-r} = -\frac{1}{r'-r}\,\int\limits_{B_{r'}\setminus B_r}\alpha_\psi\wedge\partial\rho\wedge f^\star\omega_{p-1} = -\frac{1}{r'-r}\,\int\limits_{B_{r'}\setminus B_r}\Lambda_{f^\star\omega}(\alpha_\psi\wedge\partial\rho)\, f^\star\omega_p\end{eqnarray} for all $0<r<r'$.

    Now, recalling that $\tau(z) = |z|^2$ for all $z\in\C^p$, we get $\tau = \rho^2$, hence $$i\partial\tau\wedge\bar\partial\tau = 4\rho^2\,i\partial\rho\wedge\bar\partial\rho = 4|z|^2\,i\partial\rho\wedge\bar\partial\rho,  \hspace{6ex} z\in\C^p.$$ 

Meanwhile, for any $(1,\,0)$-forms $\alpha$ and $\beta$ on a $p$-dimensional complex manifold (e.g. $\C^p$) equipped with a Hermitian metric $\gamma$, we have the general formula: \begin{eqnarray}\label{eqn:general-formula_inner-prod_trace}\langle\alpha,\,\beta\rangle_\gamma = \Lambda_\gamma(i\alpha\wedge\overline\beta),\end{eqnarray} where $\Lambda_\gamma$ is the adjoint of $\gamma\wedge\cdot$ w.r.t. the pointwise inner product $\langle\,\cdot\,,\,\cdot\,\rangle_\gamma$. This formula follows by putting together the following equalities (where $dV_\gamma:=\gamma_p$ is the volume form induced by $\gamma$): \begin{eqnarray*}\langle\alpha,\,\beta\rangle_\gamma\,dV_\gamma & = & \alpha\wedge\star_\gamma\overline\beta = i\alpha\wedge\overline\beta\wedge\gamma_{p-1}\\
\Lambda_\gamma(i\alpha\wedge\overline\beta)\,dV_\gamma & = & (i\alpha\wedge\overline\beta)\wedge\gamma_{p-1},\end{eqnarray*} where the primitivity of the $(0,\,1)$-form $\overline\beta$ yielded $\star_\gamma\overline\beta = i\overline\beta\wedge\gamma_{p-1}$ thanks to the standard formula (\ref{eqn:prim-form-star-formula-gen}) applied to the Hodge star operator $\star_\gamma$ evaluated on primitive forms.

In particular, (\ref{eqn:general-formula_inner-prod_trace}) yields the latter equality below: \begin{eqnarray}\label{eqn:d-tau_norm-trace}|d\tau|_{f^\star\omega}^2 = 2\,|\partial\tau|_{f^\star\omega}^2 = 2\,\Lambda_{f^\star\omega}(i\partial\tau\wedge\bar\partial\tau).\end{eqnarray} We also get the analogous equalities for $\rho$ in place of $\tau$.

    On the other hand, after noticing that (\ref{eqn:general-formula_inner-prod_trace}) also yields the equality $-\Lambda_{f^\star\omega}(i\alpha_\psi\wedge\partial\rho) = \langle\partial\rho,\,\overline{\alpha_\psi}\rangle_{f^\star\omega}$ (the last expression being the pointwise inner product w.r.t. $f^\star\omega$ of two $(1,\,0)$-forms), the Cauchy-Schwarz inequality yields the first of the following pointwise inequalities in $\C^p\setminus\Sigma_f$: \begin{eqnarray}\nonumber\label{eqn:H-estimate_2}\bigg|\Lambda_{f^\star\omega}(\alpha_\psi\wedge\partial\rho)\bigg| & \leq & |\alpha_\psi|_{f^\star\omega}\,|\partial\rho|_{f^\star\omega} = |\psi|_\omega\,\sqrt{\Lambda_{f^\star\omega}(i\partial\rho\wedge\bar\partial\rho)} = \frac{|\psi|_\omega}{2|z|}\,\sqrt{\Lambda_{f^\star\omega}(i\partial\tau\wedge\bar\partial\tau)} \\
      & \leq & \frac{|\partial\tau|_{f^\star\omega}}{2|z|}\,||\psi||_{C^0_\omega} = \frac{|d\tau|_{f^\star\omega}}{2\sqrt{2}\,|z|}\,||\psi||_{C^0_\omega}\end{eqnarray} for every form $\psi\in C^\infty_{p-1,\,p}(X,\,\C)$, where (\ref{eqn:f-star-psi_Hodge-star_norms}) was used to get the first equality while (\ref{eqn:d-tau_norm-trace}) and its analogue for $\partial\rho$ were also used.

    Putting together (\ref{eqn:H-estimate_1}) and (\ref{eqn:H-estimate_2}), we get: \begin{eqnarray}\label{eqn:H-estimate_3}\bigg|\frac{H(r') - H(r)}{r'-r}\bigg|\leq \frac{1}{r'-r}\,\int\limits_{B_{r'}\setminus B_r}|\Lambda_{f^\star\omega}(\alpha_\psi\wedge\partial\rho)|\, f^\star\omega_p \leq \frac{1}{r'-r}\,\frac{||\psi||_{C^0_\omega}}{2\sqrt{2}}\,\int\limits_{B_{r'}\setminus B_r}\frac{|d\tau|_{f^\star\omega}}{|z|}\, f^\star\omega_p.\end{eqnarray}

    Now, recall that $f^\star\omega_p = (d\tau/|d\tau|_{f^\star\omega})\wedge d\sigma_{\omega,\,f}$, where $d\sigma_{\omega,\,f} = \star_{f^\star\omega}(d\tau/|d\tau|_{f^\star\omega})$ is the area measure induced by the degenerate metric $f^\star\omega$ on the spheres of $\C^p$ (see (\ref{eqn:Hodge-star_degenerate}) and (\ref{eqn:area-measure_degenerate})). Thus, (\ref{eqn:H-estimate_3}) and a new application of the Fubini theorem lead to \begin{eqnarray}\label{eqn:H-estimate_4}\bigg|\frac{H(r') - H(r)}{r'-r}\bigg|\leq\frac{1}{r'-r}\,\frac{||\psi||_{C^0_\omega}}{2\sqrt{2}}\,\int\limits_r^{r'}\frac{1}{t}\,\bigg(\int\limits_{S_t}d\sigma_{\omega,\,f,\,t}\bigg)\,2tdt = \bigg(\frac{1}{r'-r}\,\int\limits_r^{r'}A_{\omega,\,f}(S_t)\,dt\bigg)\,\frac{||\psi||_{C^0_\omega}}{\sqrt{2}}\end{eqnarray} since $\tau = t^2$ on $S_t$, hence $d\tau = 2tdt$.

    Inequality (\ref{eqn:H-estimate_4}) holds for every $\psi\in C^\infty_{p-1,\,p}(X,\,\C)$ and for all reals $0<r<r'$. Fixing $r>0$ and letting $r'\downarrow r$, the left side of (\ref{eqn:H-estimate_4}) converges to $|H'(r)| = |h(r)|$, while the right side of (\ref{eqn:H-estimate_4}) converges to $A_{\omega,\,f}(S_r)\,||\psi||_{C^0_\omega}/\sqrt{2}$. Recalling (\ref{eqn:proof_main-estimate_1}), this leads to \begin{eqnarray*}|\langle\partial T_r,\,\psi\rangle| = \bigg|\bigg\langle\partial\bigg(\frac{f_\star[B_r]}{\mbox{Vol}_{\omega,\,f}(B_r)}\bigg),\,\psi\bigg\rangle\bigg|\leq\frac{1}{\sqrt{2}}\,\frac{A_{\omega,\,f}(S_r)}{\mbox{Vol}_{\omega,\,f}(B_r)}\,||\psi||_{C^0_\omega}\end{eqnarray*} for every $r>0$ and every $\psi\in C^\infty_{p-1,\,p}(X,\,\C)$. Taking the supremum over $\psi\in{\cal F}_\omega(p-1,\,p)$, we get (\ref{eqn:T_r_norm-estimate}) and we are done.  \hfill $\Box$

\vspace{3ex}

Our existence result for Ahlfors currents now follows at once from Theorem \ref{The:del-T_r-norm_estimate}.

\begin{The}\label{The:Ahlfors-current} Let $X$ be an $n$-dimensional compact complex manifold equipped with a Hermitian metric $\omega$. Suppose there exists $p\in\{1,\dots , n-1\}$ and a holomorphic map $f:\C^p\longrightarrow X$, non-degenerate at some point $x_0\in\C^p$, satisfying the condition \begin{eqnarray}\label{eqn:Ahlfors-current}\liminf\limits_{r\to +\infty}\frac{A_{\omega,\,f}(S_r)}{\mbox{Vol}_{\omega,\,f}(B_r)} = 0,\end{eqnarray} where $A_{\omega,\,f}(S_r)$ and $\mbox{Vol}_{\omega,\,f}(B_r)$ are the $(\omega,\,f)$-area of the Euclidean sphere $S_r\subset\C^p$, respectively the $(\omega,\,f)$-volume of the Euclidean ball $B_r\subset\C^p$, defined in (\ref{eqn:omega-f_area_def}) and (\ref{eqn:omega-f_vol_def}). 

Then, there exists a sequence of positive reals $r_\nu\to +\infty$ such that the currents $$T_{r_\nu}:=\frac{1}{\mbox{Vol}_{\omega,\,f}(B_{r_\nu})}\,f_\star[B_{r_\nu}]$$ converge in the weak topology of currents to a strongly positive current $T\geq 0$ of bidegree $(n-p,\,n-p)$ and of mass $1$ with respect to $\omega$ on $X$ that has the further property $dT=0$.  

\end{The}  

\noindent {\it Proof.} As explained at the beginning of this $\S.$\ref{subsection:Ahlfors_existence}, there always exists a current $T$ obtained as the weak limit of a sequence of currents $T_{r_\nu}$ as in the statement and satisfying all the stated properties except, possibly, $dT=0$. To guarantee this last property, we use hypothesis (\ref{eqn:Ahlfors-current}) to infer the existence of a sequence of positive reals $r_\nu\to +\infty$ such that \begin{eqnarray*}\label{eqn:Ahlfors-current_proof_subseq}\lim\limits_{\nu\to +\infty}\frac{A_{\omega,\,f}(S_{r_\nu})}{\mbox{Vol}_{\omega,\,f}(B_{r_\nu})} = 0.\end{eqnarray*} This implies, thanks to the general estimate (\ref{eqn:T_r_norm-estimate}), that $||\partial T_{r_\nu}||$ converges to $0$ as $\nu\to +\infty$.

Now, thanks to (\ref{eqn:unit-mass_T-r}), we have $\int_XT_{r_\nu}\wedge\omega_p = 1$ for all $\nu\in\N$. Thus, we can apply to the sequence $(T_{r_\nu})_{\nu\geq 0}$ the argument given at the beginning of this $\S.$\ref{subsection:Ahlfors_existence} to infer the existence of a subsequence (denoted by the same symbol) of $(T_{r_\nu})_{\nu\geq 0}$ that converges weakly to a current $T$. 

The only conclusion that still needs proving is $dT=0$. Since $\partial$ is continuous w.r.t. the weak topology of currents, we get: \begin{eqnarray*}0\leq||\partial T|| = ||\lim\limits_{\nu\to +\infty}\partial T_{r_\nu}||\leq\lim\limits_{\nu\to +\infty}||\partial T_{r_\nu}|| = 0,\end{eqnarray*} where the last equality was seen above. We infer that $||\partial T|| = 0$, hence $\partial T=0$, hence $dT=0$, as desired.  \hfill $\Box$

\vspace{2ex}

Since the manifold $X$ of Theorem \ref{The:Ahlfors-current} is compact, any two Hermitian metrics $\omega_1$ and $\omega_2$ thereon are comparable, in the sense that there exist constants $A, B>0$ such that $A\omega_1\leq\omega_2\leq B\omega_1$ on $X$. This implies that $Af^\star\omega_1\leq f^\star\omega_2\leq Bf^\star\omega_1$ on $\C^p$, so the growth condition (\ref{eqn:Ahlfors-current}) imposed on $f$ is independent of the choice of Hermitian metric $\omega$ on $X$.

\subsection{Examples of manifolds carrying Ahlfors currents}\label{subsection:emaples_Ahlfors} Note that hypothesis (\ref{eqn:Ahlfors-current}) is a kind of subexponential growth condition on $f$ similar, though not identical, in nature to the condition introduced in Definition \ref{Def:subexp}. In the special case where the pullback of $\omega$ under $f$ coincides with the standard K\"ahler metric $\beta$ of $\C^p$, namely \begin{eqnarray}\label{eqn:f-star-omega_beta_hypothesis}f^\star\omega = \beta = \frac{1}{2}\,\sum\limits_{j=1}^pidz_j\wedge d\bar{z}_j,\end{eqnarray} the $(\omega,\,f)$-area function $r\mapsto A_{\omega,\,f}(S_r)=A_\beta(S_r)$ for the Euclidean spheres $S_r\subset\C^p$ is the derivative of the $(\omega,\,f)$-volume function $r\mapsto \mbox{Vol}_{\omega,\,f}(B_r)=\mbox{Vol}_\beta(B_r)$ of the Euclidean balls $B_r\subset\C^p$, as shown by the Fubini theorem: \begin{eqnarray}\label{eqn:Fubini_Euclidean}\mbox{Vol}_\beta(B_r) = \int\limits_{B_r}\beta_p = \int\limits_0^r\bigg(\int\limits_{S_t}d\sigma_{\beta,\,t}\bigg)\,dt = \int\limits_0^rA_\beta(S_t)\,dt, \hspace{5ex} r>0,\end{eqnarray} where $d\sigma_{\beta,\,t}$ is the restriction to $S_t$ of the area measure $d\sigma_\beta:=\star_\beta(d\tau/|d\tau|_\beta)$ induced by the Euclideam metric $\beta$ on the spheres of $\C^p$. Thus, $\beta_p = (d\tau/|d\tau|_\beta)\wedge d\sigma_\beta$ in $\C^p$. To justify the factor $dt$ in the integrals of (\ref{eqn:Fubini_Euclidean}), note that $i\partial\tau\wedge\bar\partial\tau = \sum_{1\leq j,\,k\leq p}z_k\bar{z}_j\,idz_j\wedge d\bar{z}_k$ in $\C^p$, so (\ref{eqn:d-tau_norm-trace}) implies \begin{eqnarray*}|d\tau|^2_\beta = 2\,\Lambda_\beta(i\partial\tau\wedge\bar\partial\tau) = 4|z|^2 \hspace{6ex} \mbox{in}\hspace{1ex} \C^p.\end{eqnarray*} In particular, $|d\tau|_\beta = 2t$ on $S_t$. Meanwhile, $\tau = t^2$ on $S_t$, so $d\tau = 2t\,dt$ and $d\tau/|d\tau|_\beta = dt$ on $S_t$.

Formula (\ref{eqn:Fubini_Euclidean}) is, of course, standard, as are the equalities $\mbox{Vol}_\beta(B_r) = c_p\,r^{2p}$ and $A_\beta(S_r) = b_p\,r^{2p-1}$ with constants $c_p,b_p>0$ related by $b_p = 2p\,c_p$. In particular, hypothesis (\ref{eqn:Ahlfors-current}) is satisfied by any map $f:\C^p\longrightarrow X$ with the property (\ref{eqn:f-star-omega_beta_hypothesis}), the growth of any such $f$ being even polynomial. Thus, an immediate consequence of our Theorem \ref{The:Ahlfors-current} is the existence of an Ahlfors current induced by any map satisfying condition (\ref{eqn:f-star-omega_beta_hypothesis}).

\begin{Cor}\label{Cor:Ahlfors-current_Euclidean}Let $X$ be an $n$-dimensional compact complex manifold. Suppose there exist $p\in\{1,\dots , n-1\}$, a holomorphic map $f:\C^p\longrightarrow X$ and a Hermitian metric $\omega$ on $X$ such that $f^\star\omega = \beta$. 

Then, there exists a sequence of positive reals $r_\nu\to +\infty$ such that the currents $$T_{r_\nu}:=\frac{1}{\mbox{Vol}_\beta(B_{r_\nu})}\,f_\star[B_{r_\nu}]$$ converge in the weak topology of currents to a strongly positive current $T\geq 0$ of bidegree $(n-p,\,n-p)$ and of mass $1$ with respect to $\omega$ on $X$ that has the further property $dT=0$.  

\end{Cor} 

This corollary, in turn, implies the existence of an Ahlfors current of bidegree $(1,\,1)$ on every {\it complex torus} $X=\C^n/\Gamma$ (where $\Gamma\subset(\C^n,\,+)$ is any lattice) and on every {\it Nakamura manifold} $X=G/\Gamma$ (defined as the quotient of the solvable, non-nilpotent complex Lie group $G=(\C^3,\,\star)$ whose group operation is $$(\zeta_1,\,\zeta_2,\,\zeta_3)\star(z_1,\,z_2,\,z_3) = (\zeta_1+z_1,\,\zeta_2+e^{-\zeta_1}\,z_2,\,\zeta_3+e^{\zeta_1}\,z_3),$$ by any lattice $\Gamma\subset G$; see e.g. [Nak75] for the definition and the basic properties of these manifolds). Indeed, these manifolds $X$ were shown to be {\it non-balanced hyperbolic} in $\S2.3$(VI)(a),(c) of [MP22a] owing to the existence in each case of a non-degenerate holomorphic map $f:\C^{n-1}\longrightarrow X$ (where $n=3$ in the case of the Nakamura manifolds) satisfying property (\ref{eqn:f-star-omega_beta_hypothesis}) for a certain choice of Hermitian metric $\omega$ on each of these $X$. 

\vspace{2ex}

As a further consequence of our Theorem \ref{The:Ahlfors-current}, we now prove the existence of an Ahlfors current on the {\it Iwasawa manifold} induced by a map $f$ that does not have property (\ref{eqn:f-star-omega_beta_hypothesis}). Recall that the Iwasawa manifold is the $3$-dimensional compact complex manifold $X=G/\Gamma$ defined as the quotient of the Heisenberg group $G$, namely the nilpotent complex Lie group $G=(\C^3,\,\star)$ whose group operation is $$(\zeta_1,\,\zeta_2,\,\zeta_3)\star(z_1,\,z_2,\,z_3) = (\zeta_1+z_1,\,\zeta_2+z_2,\,\zeta_3+z_3 + \zeta_1 z_2),$$ by the lattice $\Gamma\subset G$ consisting of the elements $(z_1,\,z_2,\,z_3)\in G$ with $z_1,\,z_2,\,z_3\in\Z[i]$. (See e.g. [Nak75]).

It was observed in [MP22a, $\S2.3$(VI)(b)] that there exist a non-degenerate holomorphic map $f:\C^2\longrightarrow X$ and a Hermitian metric $\omega$ on $X$ such that \begin{eqnarray*}f^\star\omega = idz_1\wedge d\bar{z}_1 + (1+|z_1|^2)\,idz_2\wedge d\bar{z}_2  \hspace{6ex} \mbox{on}\hspace{1ex} \C^2.\end{eqnarray*} This implies that $f^\star\omega_2 = (1+|z_1|^2)\,dV_0$, where $dV_0 = idz_1\wedge d\bar{z}_1\wedge idz_2\wedge d\bar{z}_2$ is the Euclidean volume form of $\C^2$, from which we get: \begin{eqnarray*}\mbox{Vol}_{\omega,\,f}(B_r) = \int\limits_{B_r}(1+|z_1|^2)\,dV_0\leq c\,r^4(1+r^2), \hspace{5ex} r>0,\end{eqnarray*} for some constant $c>0$ independent of $r$. This implies, after fixing a constant $c_1>c$, that there exists $r_1>0$ such that \begin{eqnarray}\label{eqn:Iwasawa_vol-estimate}\mbox{Vol}_{\omega,\,f}(B_r)\leq c_1r^6, \hspace{5ex} r>r_1.\end{eqnarray}

Moreover, (\ref{eqn:d-tau_norm-trace}) gives the first equality below: \begin{eqnarray*}|d\tau|_{f^\star\omega}^2 = 2\,\Lambda_{f^\star\omega}(i\partial\tau\wedge\bar\partial\tau) = 2\,\Lambda_{f^\star\omega}\bigg(\sum\limits_{1\leq j,\,k\leq 2}z_k\bar{z}_j\,idz_j\wedge d\bar{z}_k\bigg) = 2\,\bigg(|z_1|^2 + \frac{|z_2|^2}{1+|z_1|^2}\bigg),\end{eqnarray*} for all $z=(z_1,\,z_2)\in\C^2$. When $z\in S_t$, this translates to: \begin{eqnarray}\label{eqn:Iwasawa_d-tau-estimate}|d\tau|_{f^\star\omega}^2 = 2\,\frac{|z_1|^4 + t^2}{1+|z_1|^2}\leq 2\,(|z_1|^2 + t^2)\leq 2\,(|z|^2 + t^2) = 4t^2.\end{eqnarray} 

%Thus, $|d\tau|_{f^\star\omega}\leq 2t$ on $S_t$.

For any $0<r<r'$, we get: \begin{eqnarray*}\mbox{Vol}_{\omega,\,f}(B_{r'}\setminus B_r) = \int\limits_r^{r'}\bigg(\int\limits_{S_t}\frac{1}{|d\tau|_{f^\star\omega}}\,d\sigma_{\omega,\,f,\,t}\bigg)\,2tdt \geq \int\limits_r^{r'}A_{\omega,\,f}(S_t)\,dt,\end{eqnarray*} where we used formula (\ref{eqn:comparison_measures}) (see also (\ref{eqn:Volume_Fubini})) to get the equality, while the inequality followed from $|d\tau|_{f^\star\omega}\leq 2t$ on $S_t$, a direct consequence of (\ref{eqn:Iwasawa_d-tau-estimate}). Dividing by $r'-r$ and letting $r'\downarrow r$, we further get: \begin{eqnarray*}V'(r)\geq A_{\omega,\,f}(S_r), \hspace{6ex} r>0,\end{eqnarray*} where $V'$ is the derivative of the function $r\mapsto V(r):= \mbox{Vol}_{\omega,\,f}(B_r)$. Thus, we get: \begin{eqnarray}\label{eqn:Iwasawa_A-V_ratio}\frac{A_{\omega,\,f}(S_r)}{\mbox{Vol}_{\omega,\,f}(B_r)}\leq\frac{V'(r)}{V(r)}, \hspace{6ex} r>0.\end{eqnarray}

Now, fix $\alpha>6$. We claim that there exists a sequence $r_\nu\longrightarrow +\infty$ such that \begin{eqnarray}\label{eqn:Iwasawa_proof-claim}\frac{V'(r_\nu)}{V(r_\nu)}<\frac{\alpha}{r_\nu},   \hspace{6ex} \nu\in\N.\end{eqnarray} We prove this claim by contradiction. Suppose the contention is false. Then, there exists $r_2>0$ such that $V'(r)/V(r)\geq\alpha/r$ for all $r\geq r_2$. This is equivalent to $(\log V(r))'\geq (\alpha\,\log r)'$, hence to $(\log(V(r)/r^\alpha))'\geq 0$ for all $r\geq r_2$. This implies that the function $[r_2,\,+\infty)\ni r\mapsto\log(V(r)/r^\alpha)$ is non-decreasing, hence $\log(V(r)/r^\alpha)\geq\log(V(r_2)/r_2^\alpha):=c_2$ for all $r\geq r_2$. We would then have $V(r)\geq C_2\,r^\alpha$ for all $r\geq r_2$, where $C_2:=\exp(c_2)$. However, this would contradict (\ref{eqn:Iwasawa_vol-estimate}).

Thus, the claim is proved. Putting (\ref{eqn:Iwasawa_A-V_ratio}) and (\ref{eqn:Iwasawa_proof-claim}) together, we get: \begin{eqnarray*}\frac{A_{\omega,\,f}(S_{r_\nu})}{\mbox{Vol}_{\omega,\,f}(B_{r_\nu})}\leq\frac{V'(r_\nu)}{V(r_\nu)}<\frac{\alpha}{r_\nu}\longrightarrow 0 \hspace{3ex}\mbox{as}\hspace{1ex}\nu\longrightarrow +\infty.\end{eqnarray*} This proves that \begin{eqnarray*}\liminf\limits_{r\to +\infty}\frac{A_{\omega,\,f}(S_r)}{\mbox{Vol}_{\omega,\,f}(B_r)} = 0,\end{eqnarray*} which means that the non-degenerate holomorphic map $f:\C^2\longrightarrow X$ satisfies condition (\ref{eqn:Ahlfors-current}). 

Thus, our Theorem \ref{The:Ahlfors-current} yields the following 

\begin{Prop}\label{Prop:Iwasawa_Ahlfors} The Iwasawa manifold carries an Ahlfors current.

\end{Prop}  

\subsection{Link with partial hyperbolicity}\label{subsection:link_hyp-Ahlfors} When $p=n-1$, a holomorphic map $f:\C^{n-1}\longrightarrow (X,\,\omega)$ that is non-degenerate at some point and satisfies the growth condition (\ref{eqn:Ahlfors-current}) gives rise, thanks to Theorem \ref{The:Ahlfors-current}, to an Ahlfors current $T\geq 0$ of bidegree $(1,\,1)$ on $X$. Since $T$ is a non-zero (because it has unit mass w.r.t. $\omega$) $d$-closed positive $(1,\,1)$-current, $X$ cannot carry any {\it degenerate balanced structure}, namely there is no $C^\infty$ positive definite $(n-1,\,n-1)$-form $\Omega$ on $X$ such that $\Omega\in\mbox{Im}\,d$. This is one implication of an equivalence proved in [Pop15, Proposition 5.4]. Since the existence of a degenerate balanced structure on $X$ is a special type of balanced hyperbolicity (cf. [MP22a]), this observation is an indication (already noted e.g. in [MQ98] in the context of maps from $\C$) of a link between possible hyperbolicity properties of $X$ and the possible non-existence of Ahlfors currents thereon.

We will now make this link precise for an arbitrary $p\in\{1,\dots , n-1\}$. Our result in this direction will be similar in nature to Theorem \ref{The:partial-p-hyperbolicity_implication}, except that the subexponential growth condition of Definition \ref{Def:subexp} will be replaced by the growth condition (\ref{eqn:Ahlfors-current}). This leads naturally to another notion of partial hyperbolicity.

  \begin{Def}\label{Def:strongly_partial-p-hyperbolic} Let $X$ be a compact complex manifold with $\mbox{dim}_\C X =n$. Fix $p\in\{1,\dots , n-1\}$ and suppose there is a $C^\infty$ complex vector subbundle $E$ of rank $\geq p$ of the holomorphic tangent bundle $T^{1,\,0}X$.

    The manifold $X$ is said to be {\bf strongly partially $p$-hyperbolic} in the $E$-directions if there exists no holomorphic map $f:\C^p\longrightarrow X$ simultaneously satisfying the following three conditions:

\vspace{1ex}

(i)\, $f$ is non-degenerate at some point $x_0\in\C^p$;

\vspace{1ex}

(ii)\, $f$ is $E$-horizontal;

\vspace{1ex}

    (iii)\, for some (hence any) Hermitian metric $\omega$ on $X$, $f$ satisfies the growth condition (\ref{eqn:Ahlfors-current}), namely \begin{eqnarray*}\liminf\limits_{r\to +\infty}\frac{A_{\omega,\,f}(S_r)}{\mbox{Vol}_{\omega,\,f}(B_r)} = 0,\end{eqnarray*} where $A_{\omega,\,f}(S_r)$ and $\mbox{Vol}_{\omega,\,f}(B_r)$ are the $(\omega,\,f)$-area of the Euclidean sphere $S_r\subset\C^p$, respectively the $(\omega,\,f)$-volume of the Euclidean ball $B_r\subset\C^p$, defined in (\ref{eqn:omega-f_area_def}) and (\ref{eqn:omega-f_vol_def}).

  \end{Def}

  \vspace{2ex}

  We can now state the result we have been alluding to.

\begin{The}\label{The:partial-p-K-hyperbolicity-new-growth_implication} Let $X$ be an $n$-dimensional compact complex manifold and let $p\in\{1,\dots , n-1\}$.

  If $X$ has a {\bf partially $p$-K\"ahler hyperbolic structure} $(E,\,\Omega,\,\omega = \omega_E + \omega_{nE})$, then $X$ is {\bf strongly partially $p$-hyperbolic} in the $E$-directions. In particular, there is no Ahlfors current on $X$ induced by a map as in Definition \ref{Def:strongly_partial-p-hyperbolic}.

\end{The}

\noindent {\it Proof.} We reason by contradiction. Suppose that an $E$-horizontal holomorphic map $f:\C^p\longrightarrow X$ that is non-degenerate at some point and satisfies the growth condition (\ref{eqn:Ahlfors-current}) existed. We use the notation of $\S.$\ref{section:partial-hyp_two-definitions}. 

Since $f$ is $E$-horizontal, $f^\star\omega = f^\star\omega_E$ (cf. (\ref{eqn:f-star_omega_E})), hence also $f^\star\omega_p = f^\star\Omega = d(\widetilde{f}^\star\widetilde\Gamma)$ on $\C^p$, where $\widetilde\Gamma$ is the $\widetilde\omega$-bounded $(2p-1)$-form on the universal cover $\widetilde{X}$ with the property $\pi_X^\star\Omega = d\widetilde\Gamma$ given by the partial $p$-K\"ahler hyperbolicity hypothesis on $X$. Thus, we get: \begin{eqnarray*}1 & = & \frac{1}{\mbox{Vol}_{\omega,\,f}(B_r)}\,\int\limits_{B_r}f^\star\omega_p = \frac{1}{\mbox{Vol}_{\omega,\,f}(B_r)}\,\int\limits_{B_r}d(\widetilde{f}^\star\widetilde\Gamma) = \frac{1}{\mbox{Vol}_{\omega,\,f}(B_r)}\,\int\limits_{S_r}\widetilde{f}^\star\widetilde\Gamma \\
& \leq & \frac{C}{\mbox{Vol}_{\omega,\,f}(B_r)}\,\int\limits_{S_r}\sigma_{\omega,\,f,\,r} = C\,\frac{A_{\omega,\,f}(S_r)}{\mbox{Vol}_{\omega,\,f}(B_r)},  \hspace{6ex} r>0,\end{eqnarray*} where $C>0$ is a constant whose existence reflects the boundedness in $\C^p$ of $\widetilde{f}^\star\widetilde\Gamma$ w.r.t. $f^\star\omega$ (cf. Lemma \ref{Lem:Gamma-tilde_bounded}).

This contradicts the assumption (\ref{eqn:Ahlfors-current}) made on $f$. \hfill $\Box$

\vspace{2ex}

  The use of the adverb ``strongly'' in the above terminology is warranted by the partial hyperbolicity notion introduced in Definition \ref{Def:strongly_partial-p-hyperbolic} being stronger than the one introduced in Definition \ref{Def:partial_p-hyperbolic}, as shown by the following

\begin{The}\label{The:strongly-partially-hyp_implies_partially-hyp} Let $X$ be a compact complex manifold with $\mbox{dim}_\C X =n\geq 3$. Fix $p\in\{1,\dots , n-1\}$ and suppose there exists a $C^\infty$ complex vector subbundle $E$ of rank $\geq p$ of the holomorphic tangent bundle $T^{1,\,0}X$.

  If $X$ is {\bf strongly partially $p$-hyperbolic} in the $E$-directions, then $X$ is {\bf partially $p$-hyperbolic} in the $E$-directions.

\end{The}

The proof of this result follows at once from the following general fact.

\begin{Prop}\label{Prop:growth-notions_implication} Let $X$ be an $n$-dimensional compact complex manifold, let $p\in\{1,\dots , n-1\}$ and let $f:\C^p\longrightarrow X$ be a holomorphic map that is non-degenerate at some point $x_0\in\C^p$.

  If $f$ satisfies the subexponential growth condition of Definition \ref{Def:subexp}, then $f$ satisfies the growth condition (\ref{eqn:Ahlfors-current}).

\end{Prop}  

\noindent {\it Proof.} This proof is implicit in the proof of Theorem \ref{The:partial-p-hyperbolicity_implication}. Indeed, fix an arbitrary Hermitian metric $\omega$ on $X$. As already pointed out, inequality (\ref{eqn:Volume_1}) holds without any special assumption on $f$ or $X$.

On the other hand, let us suppose that $f$ does not satisfy the growth condition (\ref{eqn:Ahlfors-current}). This means that there exists a constant $C>0$ such that \begin{eqnarray*}\frac{A_{\omega,\,f}(S_r)}{\mbox{Vol}_{\omega,\,f}(B_r)}\geq\frac{1}{C}\end{eqnarray*} for every $r>0$ large enough. This amounts to inequality (\ref{eqn:Volume_Stokes}) (which in the proof of Theorem \ref{The:partial-p-hyperbolicity_implication} is a consequence of the Stokes theorem and the partial p-K\"ahler hyperbolicity assumption on $X$) holding for all $r>0$.

As in the proof of Theorem \ref{The:partial-p-hyperbolicity_implication}, putting (\ref{eqn:Volume_1}) and (\ref{eqn:Volume_Stokes}) together and using part (i) of the subexponential growth assumption on f, we get inequality (\ref{eqn:Volume_final}) which, thanks to part (ii) of the subexponential growth assumption on f, leads to a contradiction of the non-degeneracy property of $f$.  \hfill $\Box$

%\vspace{3ex}

%\noindent {\it Proof of Theorem \ref{The:strongly-partially-hyp_implies_partially-hyp}.} Proving the contention is equivalent to proving the implication:

%\vspace{1ex}

%{\it If $X$ is not partially $p$-hyperbolic in the $E$-directions, then $X$ is not strongly partially $p$-hyperbolic in the $E$-directions.}

%\vspace{1ex}

%This implication is an immediate consequence of Proposition \ref{Prop:growth-notions_implication} and of Definitions \ref{Def:strongly_partial-p-hyperbolic} and \ref{Def:partial_p-hyperbolic}.  \hfill $\Box$

\section{Curvature sign and partial hyperbolicity}\label{section:curvature} The curvature-like notion that seems best suited to our partial hyperbolicity context and will be used in this paper was introduced very recently in [Pop23]. We now briefly recall the construction for the reader's convenience.

\subsection{Brief reminder of a construction from [Pop23]}\label{subsection:reminder_pluriclosed-star-split}

If $\omega$ is a Hermitian metric on an $n$-dimensional complex manifold $X$, the pointwise linear map \begin{eqnarray*}\omega_{n-2}\wedge\cdot :\Lambda^{1,\,1}T^\star X\longrightarrow\Lambda^{n-1,\,n-1}T^\star X\end{eqnarray*} (that multiplies $(1,\,1)$-forms by $\omega_{n-2}$) is {\it bijective}. Therefore, there exists a unique $C^\infty$ $(1,\,1)$-form $\rho_\omega$ on $X$ such that \begin{eqnarray*}\label{eqn:rho_omega_def}i\partial\bar\partial\omega_{n-2} = \omega_{n-2}\wedge\rho_\omega.\end{eqnarray*}

  Explicit computations (cf. [Pop23, $\S2.1$]) lead to the following formula for the Hodge star operator $\star=\star_\omega$ evaluated on $\rho_\omega$: \begin{eqnarray*}\label{eqn:star-rho_omega_formula}\star\rho_\omega = \frac{1}{n-1}\,\frac{\omega\wedge i\partial\bar\partial\omega_{n-2}}{\omega_n}\,\omega_{n-1} - i\partial\bar\partial\omega_{n-2}.\end{eqnarray*} This prompts one to consider the $C^\infty$ function $f_\omega:X\to\R$ defined by \begin{eqnarray}\label{eqn:function_f_omega}f_\omega:=\frac{\omega\wedge i\partial\bar\partial\omega_{n-2}}{\omega_n}.\end{eqnarray} It is uniquely associated with each Hermitian metric $\omega$ on $X$ and plays a role somewhat similar in some respects, though quite different in other respects, to that of the classical scalar curvature. Several results of [Pop23] suggest that the scalar-valued $(1,\,1)$-form $\rho_\omega$ and the scalar-valued $(n-1,\,n-1)$-form \begin{eqnarray}\label{eqn:star-rho_omega_formula}\star\rho_\omega = \frac{1}{n-1}\,f_\omega\,\omega_{n-1} - i\partial\bar\partial\omega_{n-2}\end{eqnarray} play roles analogous in some respects to those of the classical vector bundle-valued curvature form $i\Theta_\omega(T^{1,\,0}X)\in C^\infty_{1,\,1}(X,\,\mbox{End}\,(T^{1,\,0}X))$ of the holomorphic tangent bundle of $X$.

\vspace{1ex}

  We also recall the following definition from [Pop23]:

  \vspace{1ex}

  {\it The Hermitian metric $\omega$ on $X$ is said to be {\bf pluriclosed star split} if $\partial\bar\partial(\star\rho_\omega) =0$.}

  \vspace{1ex}

  The following statement reproduces several results from [Pop23] that will be relevant to the ensuing discussion despite them not being used directly therein.

  \begin{The}\label{The:pluriclosed-star-split_summary} Let $X$ be a compact connected complex manifold with $\mbox{dim}_\C X = n\geq 3$ and let $\omega$ be a Hermitian metric on it.

    \vspace{1ex}

 (i) ([Pop23, Proposition 1.1.])\, If $\omega$ is {\bf pluriclosed star split}, its associated function $f_\omega$ satisfies one of the following three conditions: \begin{eqnarray*}f_\omega> 0 \hspace{1ex} \mbox{on} \hspace{1ex} X \hspace{5ex} \mbox{or} \hspace{5ex} f_\omega< 0 \hspace{1ex} \mbox{on} \hspace{1ex} X  \hspace{5ex} \mbox{or} \hspace{5ex} f_\omega\equiv 0.\end{eqnarray*}    

\vspace{1ex}

(ii) ([Pop23, Proposition 2.5.])\, If the function $f_\omega$ is a {\bf non-zero constant}, the metric $\omega$ is {\bf pluriclosed star split} if and only if it is {\bf Gauduchon}.

\vspace{1ex}

(iii) ([Pop23, Proposition 2.5.])\, If the metric $\omega$ is {\bf Gauduchon}, then it is {\bf pluriclosed star split} if and only if the function $f_\omega$ is {\bf constant}.

\vspace{1ex}

(iv) ([Pop23, Proposition 1.2.])\, If $\omega$ is {\bf balanced} and {\bf pluriclosed star split}, the (necessarily constant) function $f_\omega$ is {\bf non-negative}.

  \end{The}

  \vspace{1ex}

  Recall that a Hermitian metric $\omega$ on an $n$-dimensional complex manifold $X$ is said to be {\it Gauduchon} (see [Gau77a]), respectively {\it balanced} (see [Gau77b] where these metrics were termed {\it semi-K\"ahler}), if $\partial\bar\partial\omega^{n-1} = 0$, respectively if $d\omega^{n-1} = 0$. Every balanced metric is, obviously, Gauduchon. Meanwhile, Gauduchon metrics always exist on every compact complex manifold (see [Gau77a]), while balanced metrics need not exist, but they exist on quite a few classes of compact complex non-K\"ahler manifolds. Meanwhile, every K\"ahler metric is, obviously, balanced.

  \vspace{1ex}

  As evidenced by the above Theorem \ref{The:pluriclosed-star-split_summary}, the sign of the function $f_\omega$ associated with a given Hermitian metric $\omega$ gives geometric information on the manifold $X$. A question that was raised in [Pop23] and whose answer is affirmative on nilmanifolds of complex dimension $3$ is whether (or when) the sign of $f_\omega$ depends only on the complex structure of $X$.

  \vspace{1ex}

  In this paper, we shall use the signs of the function $f_\omega$ and the $(n-1,\,n-1)$-form $\star\rho_\omega$ to give a sufficient criterion for partial hyperbolicity of an $n$-dimensional compact complex Hermitian manifold $(X,\,\omega)$.

  \subsection{Sufficient metric criterion for partial hyperbolicity}\label{subsection:sufficient_metric-hyp} In this subsection, we discuss partial hyperbolicity in the directions of a vector subbundle $E\subset T^{1,\,0}X$ of complex co-rank $1$. However, we start with a definition for an arbitrary rank. 

  \begin{Def}\label{Def:partial-negativity} Let $X$ be a compact complex manifold with $\mbox{dim}_\C X=n$. Fix any $p\in\{1,\dots , n-1\}$ and suppose there is a $C^\infty$ complex vector subbundle $E\subset T^{1,\,0}X$ of rank $p$. 

A real $(p,\,p)$-form $\Omega$ on $X$ is said to be {\bf negative in the $E$-directions} if for every point $x_0\in X$ and for some (hence any) $C^\infty$ frame $\xi_1,\dots ,\xi_p$ of $E$ in a neighbourhood $U$ of $x_0$, we have $$\Omega(\xi_1,\bar\xi_1,\dots , \xi_p,\bar\xi_p)<0$$ at every point $x\in U$. 

\end{Def}    

\vspace{1ex}

  We are now in a position to state our sufficient criterion for partial hyperbolicity in terms of the curvature-like objects $f_\omega$ and $\star\rho_\omega$ associated with a Hermitian metric $\omega$ via formulae (\ref{eqn:function_f_omega}) and (\ref{eqn:star-rho_omega_formula}).

  \begin{The}\label{The:sufficient_metric-hyp} Let $X$ be an $n$-dimensional compact complex manifold such that there is a $C^\infty$ complex vector subbundle $E\subset T^{1,\,0}X$ of rank $n-1$. 

    If there exists a Hermitian metric $\omega$ on $X$ such that $f_\omega>0$ and the $(n-1,\,n-1)$-form $\star\rho_\omega$ is {\bf negative in the $E$-directions}, $X$ is {\bf strongly partially $(n-1)$-hyperbolic} in the E-directions.

\end{The}

  \noindent {\it Proof.} We reason by contradiction. Suppose there exists an $E$-horizontal holomorphic map $f:\C^{n-1}\longrightarrow X$ that is non-degenerate at some point and satisfies the growth condition (\ref{eqn:Ahlfors-current}). Then, by Theorem \ref{The:Ahlfors-current}, $f$ induces an Ahlfors current $T$ on $X$ which is the weak limit of a sequence $(T_{r_\nu})_{\nu\in\N}$ of currents of bidegree $(1,\,1)$ on $X$ defined by (\ref{eqn:T_r_def}).

  Thus, on the one hand, the property $dT=0$ and formula (\ref{eqn:star-rho_omega_formula}) imply, via the Stokes theorem, the equality below: \begin{eqnarray}\label{eqn:integral_Ahlfors-current_star-rho_1}\int\limits_X T\wedge\star\rho_\omega = \frac{1}{n-1}\,\int\limits_X f_\omega\,T\wedge\omega_{n-1} >0,\end{eqnarray} where the inequality follows from the hypothesis $f_\omega>0$ and the property $T\geq 0$ (with $T\neq 0$).

  On the other hand, since $f$ is $E$-horizontal and $\star\rho_\omega$ is negative in the $E$-directions, $f^\star(\star\rho_\omega)$ is a negative $(n-1,\,n-1)$-form on $\C^{n-1}$. We get \begin{eqnarray}\label{eqn:integral_Ahlfors-current_star-rho_2}\int\limits_X T\wedge\star\rho_\omega =\lim\limits_{\nu\to +\infty} \int\limits_X T_{r_\nu}\wedge\star\rho_\omega = \lim\limits_{\nu\to +\infty}\frac{1}{\mbox{Vol}_{\omega,\,f}(B_{r_\nu})}\,\int\limits_{B_{r_\nu}}f^\star(\star\rho_\omega) \leq 0.\end{eqnarray}

 Since (\ref{eqn:integral_Ahlfors-current_star-rho_1}) and (\ref{eqn:integral_Ahlfors-current_star-rho_2}) contradict each other, we are done. \hfill $\Box$ 

  \vspace{3ex}

  We now give applications of Theorem \ref{The:sufficient_metric-hyp} to several classes of examples of compact complex manifolds.

\vspace{1ex}

{\it (I)\, The Iwasawa manifold}

\vspace{1ex}

We now return briefly to this well-known $3$-dimensional compact complex manifold $X$ whose description was sketched in $\S$\ref{subsection:emaples_Ahlfors}. Its cohomology is completely determined by three holomorphic $(1,\,0)$-forms $\alpha$, $\beta$, $\gamma$ (so, $\bar\partial\alpha = \bar\partial\beta = \bar\partial\gamma = 0$) induced respectively by the forms $dz_1$, $dz_2$ and $dz_3 - z_1 dz_2$ on $\C^3$ by passage to the quotient. They satisfy the structure equations: \begin{eqnarray*}\label{eqn:Iwasawa_structure-eq}\partial\alpha = \partial\beta = 0  \hspace{3ex}\mbox{and} \hspace{3ex} \partial\gamma = -\alpha\wedge\beta\neq 0.\end{eqnarray*}

Let $\{\xi,\eta,\nu\}$ be the global holomorphic frame of $T^{1,\,0}X$ that is dual to the global holomorphic frame $\{\alpha,\beta,\gamma\}$ of $\Lambda^{1,\,0}T^\star X$. We consider the following holomorphic vector subbundle of $T^{1,\,0}X$: $$E:=\langle\xi,\eta\rangle.$$

The standard Hermitian metric on the Iwasawa manifold is \begin{eqnarray*}\omega = i\alpha\wedge\bar\alpha + i\beta\wedge\bar\beta + i\gamma\wedge\bar\gamma.\end{eqnarray*} Moreover, by $\S2.2.1$ of [Pop23], for this metric one has: $$f_\omega = 1 \hspace{2ex}\mbox{and}\hspace{2ex} \star\rho_\omega = \frac{1}{2}\,f_\omega\,\omega_2 - i\partial\bar\partial\omega = \frac{1}{2}\,\bigg(i\alpha\wedge\bar\alpha\wedge i\gamma\wedge\bar\gamma + i\beta\wedge\bar\beta\wedge i\gamma\wedge\bar\gamma - i\alpha\wedge\bar\alpha\wedge i\beta\wedge\bar\beta\bigg).$$ So, in particular, $\star\rho_\omega$ is {\it negative in the $E$-directions}.

  Thus, the hypotheses of Theorem \ref{The:sufficient_metric-hyp} are satisfied. Consequently, we conclude the following

\begin{Cor}\label{Cor:Iwasawa_partial-hyp} The Iwasawa manifold is {\bf strongly partially $2$-hyperbolic} in the $E$-directions.

\end{Cor}  

Recall that, on the other hand, the Iwasawa manifold $X=G/\Gamma$ was shown in $\S2.3.(VI)(b)$ of [MP22a] to not be {\it divisorially hyperbolic}. The non-degenerate holomorphic map $f:\C^2\longrightarrow X$ that was used to prove this statement in [MP22a] (and again in the proof of the above Proposition  \ref{Prop:Iwasawa_Ahlfors}) was the composition of the map $j:\C^2\longrightarrow\C^3$ defined by $j(z_1,z_2)= (z_1, z_2, 0)$ with the projection map $\pi:\C^3=G\longrightarrow X=G/\Gamma$. We wish to stress that $f$ is not $E$-horizontal, so there is no clash with Corollary \ref{Cor:Iwasawa_partial-hyp}. 

Indeed, since $$\pi^\star\alpha = dz_1,\hspace{3ex} \pi^\star\beta = dz_2,\hspace{3ex} \pi^\star\gamma = dz_3 - z_1\,dz_2,$$ we get: \begin{eqnarray*}\pi^\star(\star\rho_\omega) & = & \frac{1}{2}\,\bigg(idz_1\wedge d\bar{z}_1\wedge i(dz_3 - z_1\,dz_2)\wedge(d\bar{z}_3 - \bar{z}_1\,d\bar{z}_2) \\ 
& + & idz_2\wedge d\bar{z}_2\wedge i(dz_3 - z_1\,dz_2)\wedge(d\bar{z}_3 - \bar{z}_1\,d\bar{z}_2) -idz_1\wedge d\bar{z}_1\wedge idz_2\wedge d\bar{z}_2\bigg).\end{eqnarray*} Moreover, since the third component of $j$ vanishes, $j^\star(dz_3) = 0$ and $j^\star(d\bar{z}_3) = 0$. Thus, we get: \begin{eqnarray*}f^\star(\star\rho_\omega) = j^\star\bigg(\pi^\star(\star\rho_\omega)\bigg) = \frac{1}{2}\,(|z_1|^2 - 1)\, idz_1\wedge d\bar{z}_1\wedge idz_2\wedge d\bar{z}_2\end{eqnarray*} in $\C^2$. This shows that the $(2,\,2)$-form $f^\star(\star\rho_\omega)$ is not negative at every point of $\C^2$ as would be the case if the map $f$ were $E$-horizontal.

\vspace{2ex}

{\it (II)\, Certain small deformations $X_t$ of the Iwasawa manifold $X=X_0$}

\vspace{1ex}

According to $\S2.2.3$ of [Pop23], for the standard Hermitian metric $\omega_t$ defined by the three standard smooth $(1,\,0)$-forms $\alpha_t$, $\beta_t$, $\gamma_t$ that generate the cohomology of any small deformation $X_t$ lying in one of Nakamura's classes $(ii)$ or $(iii)$ of the Iwasawa manifold $X=X_0$, one gets: \begin{eqnarray*}f_{\omega_t} & = & A(t), \\
  \star\rho_{\omega_t} & = & \frac{A(t)}{2}\,\bigg(i\alpha_t\wedge\bar\alpha_t\wedge i\gamma_t\wedge\bar\gamma_t + i\beta_t\wedge\bar\beta_t\wedge i\gamma_t\wedge\bar\gamma_t - i\alpha_t\wedge\bar\alpha_t\wedge i\beta_t\wedge\bar\beta_t\bigg),\end{eqnarray*} where \begin{eqnarray*}A(t):=|\sigma_{12}(t)|^2 + |\sigma_{2\bar{1}}(t)|^2 + |\sigma_{1\bar{2}}(t)|^2 - 2\,\mbox{Re}\,\bigg(\sigma_{1\bar{1}}(t)\,\overline{\sigma_{2\bar{2}}(t)}\bigg)\end{eqnarray*} is a constant (depending only on the deformation parameter $t$) defined by the coefficients $\sigma_{12}(t)$ and $\sigma_{j\bar{k}}(t)$ featuring in the structure equations satisfied by the forms $\alpha_t$, $\beta_t$, $\gamma_t$.

We now see that the hypotheses of the above Theorem \ref{The:sufficient_metric-hyp}, when the vector bundle $E_t$ is chosen analogously to that of the Iwasawa manifold case (I) discussed above, are satisfied whenever $A(t)>0$. We conclude the following

\begin{Cor}\label{Cor:Iwasawa-def_partial-hyp} For every complex number $t$ close enough to $0$ and such that $A(t)>0$, the small deformation $X_t$ lying in one of Nakamura's classes (ii) or (iii) of the Iwasawa manifold $X=X_0$ is {\bf strongly partially $2$-hyperbolic} in the $E_t$-directions.
  
\end{Cor}

\vspace{2ex}

{\it (III)\, Certain small deformations $X_t$ of the Calabi-Eckmann manifold $X=(S^3\times S^3,\,J_{CE})$}

\vspace{1ex}

According to $\S2.2.5$ of [Pop23], for the standard Hermitian metric $\omega_t$ on the small deformation $X_t$ of the manifold $X=X_0$ (defined to be $S^3\times S^3$ equipped with the Calabi-Eckmann complex structure $J_{CE}$, where $S^3$ is the $3$-sphere), we have: $$f_{\omega_t} = 8\,\mbox{Im}\,(t) \hspace{2ex}\mbox{and}\hspace{2ex} \star\rho_{\omega_t} = \mbox{Im}\,(t)\,\bigg(\widehat{i\varphi_t^1\wedge\overline\varphi_t^1} + \widehat{i\varphi_t^2\wedge\overline\varphi_t^2} - \widehat{i\varphi_t^3\wedge\overline\varphi_t^3}\bigg)$$ for every $t$ in a small enough neighbourhood of $0$ in $\C$. This shows that, for every such $t$ for which $\mbox{Im}\,(t)>0$, we have: $f_{\omega_t}>0$ and $\star\rho_{\omega_t}$ is {\it negative in the $E_t$-directions}, where $E_t\subset T^{1,\,0}X_t$ is the $C^\infty$ subbundle of $T^{1,\,0}X_t$ generated by $\theta^1_t$ and $\theta^2_t$, the first two $(1,\,0)$-vector fields in the frame $\{\theta^1_t,\,\theta^2_t,\,\theta^3_t\}$ of $T^{1,\,0}X_t$ that is dual to $\{\varphi_t^1,\,\varphi_t^2,\,\varphi_t^3\}$.

Applying again Theorem \ref{The:sufficient_metric-hyp}, we conclude the following

\begin{Cor}\label{Cor:C-E_def_partial-hyp} For every complex number $t$ close enough to $0$ and such that $\mbox{Im}\,(t)>0$, the small deformation $X_t$ of the Calabi-Eckmann manifold $X=X_0= (S^3\times S^3,\,J_{CE})$ is {\bf strongly partially $2$-hyperbolic} in the $E_t$-directions.
  
\end{Cor}

\vspace{1ex}

%\newpage

\vspace{3ex}

\noindent {\bf References.} \\

\noindent [AB91]\, L. Alessandrini, G. Bassanelli --- {\it Compact $p$-K\"ahler Manifolds} --- Geom. Dedicata {\bf 38} (1991), 199–210.

\vspace{1ex}

\noindent [ADOS22]\, D. Angella, A. Dubickas, A. Otiman, J. Stelzig --- {\it On Metric and Cohomological Properties of Oeljeklaus-Toma Manifolds} --- arXiv e-print DG 2201.06377 to appear in Publicacions Matematiques.

\vspace{1ex}

\noindent [Bro78]\, R. Brody --- {\it Compact Manifolds and Hyperbolicity} --- Trans. Amer. Math. Soc. {\bf 235} (1978), 213-219.

\vspace{1ex}

\noindent [BS12]\, D. Burns, N. Sibony --- {\it Limit Currents and Value Distribution of Holomorphic Maps} --- Ann. Inst. Fourier, Grenoble, {\bf 62}, 1 (2012), 145-176.

%\vspace{1ex}

%\noindent [CY17]\, B.-L. Chen, X. Yang --- {\it Compact K\"ahler Manifolds Homotopic to Negatively Curved Riemannian Manifolds} --- Math. Ann. DOI 10.1007/s00208-017-1521-7

\vspace{1ex}

\noindent [Dem97]\, J.-P. Demailly --- {\it Complex Analytic and Algebraic Geometry} --- http://www-fourier.ujf-grenoble.fr/~demailly/books.html

\vspace{1ex}

\noindent [Dem21]\, J.-P. Demailly --- {\it A Simple Proof of the Kobayashi Conjecture on the Hyperbolicity of General Algebraic Hypersurfaces} --- Panoramas et Synth\`eses, {\bf 56} (2021), 89-134.

\vspace{1ex}

\noindent [DO98]\, S. Dragomir, L. Ornea  --- {\it Locally Conformal K\"ahler Geometry} ---  Progress in Mathematics, {\bf 155}. Birkh\"auser Boston, Inc., Boston, MA, 1998.

\vspace{1ex}

\noindent [dT10]\, H. de Th\'elin --- {\it Ahlfors' Currents in Higher Dimension} --- Ann. Fac. Sci. Toulouse Math. (6) {\bf 19} (2010), no. 1, 121-133. 

\vspace{1ex}

\noindent [FKV15]\, A. Fino, H. Kasuya, L. Vezzoni, --- {\it SKT and Tamed Symplectic Structures on Solvmanifolds} --- Tohoku Math. J. (2), Vol. {\bf 67}, No. 1 (2015), 19-37.

\vspace{1ex}

\noindent [Gau77a]\, P. Gauduchon --- {\it Le th\'eor\`eme de l'excentricit\'e nulle} --- C. R. Acad. Sci. Paris, S\'er. A, {\bf 285} (1977), 387-390.

\vspace{1ex}

\noindent [Gau77b]\, P. Gauduchon --- {\it Fibr\'es hermitiens \`a endomorphisme de Ricci non n\'egatif} --- Bull. Soc. Math. France {\bf 105} (1977) 113-140.

\vspace{1ex}

\noindent [Gro91]\, M. Gromov --- {\it K\"ahler Hyperbolicity and $L^2$ Hodge Theory} --- J. Diff. Geom. {\bf 33} (1991), 263-292.

\vspace{1ex}

\noindent [KO05] Y. Kamishima,  L. Ornea---{\it Geometric Flow on Compact Locally Conformally K\"ahler Manifolds} --- Tohoku Math. J. (2) {\bf 57} (2005), no. 2, 201-221.

\vspace{1ex}

\noindent [Kas13] H. Kasuya ---{\it Vaisman Metrics on Solvmanifolds and Oeljeklaus-Toma Manifolds} --- Bull. Lond. Math. Soc. {\bf 45} (2013), no. 1, 15-26.

\vspace{1ex}

\noindent [Kas17] H. Kasuya ---{\it Generalized Deformations and Holomorphic Poisson Cohomology of Solvmanifolds}--- Ann. Global Anal. Geom. {\bf 51} (2017), no. 2, 155-177

\vspace{1ex}

\noindent [Kas21] H. Kasuya ---{\it Remarks on Dolbeault Cohomology of Oeljeklaus-Toma Manifolds and Hodge Theory} ---  Proc. Amer. Math. Soc. {\bf 149} (2021), no. 7, 3129-3137.

\vspace{1ex}

\noindent [Kas23] H. Kasuya --- {\it Cohomology of Holomorphic Line Bundles and Hodge Symmetry on Oeljeklaus-Toma Manifolds}---   Eur. J. Math. {\bf 9} (2023), no. 3, Paper No. 50, 15 pp..

\vspace{1ex}

\noindent [Kob67]\, S. Kobayashi --- {\it Invariant Distances on Complex Manifolds and Holomorphic Mappings} --- J. Math. Soc. Japan, {\bf 19}, No. 4 (1967), 460-480.

\vspace{1ex}

\noindent [Kob70]\, S. Kobayashi --- {\it Hyperbolic Manifolds and Holomorphic Mappings} --- Marcel Dekker, New York (1970). 

%\vspace{1ex}

%\noindent [MP21]\, S. Marouani, D. Popovici --- {\it Some Properties of Balanced Hyperbolic Compact Complex Manifolds} --- Internat. J. Math., {\bf 33}, No. 3 (2022) 2250019, DOI : 10.1142/S0129167X22500197.

\vspace{1ex}

\noindent [MP22a]\, S. Marouani, D. Popovici --- {\it Balanced Hyperbolic and Divisorially Hyperbolic Compact Complex Manifolds} --- arXiv e-print CV 2107.08972v2, to appear in Mathematical Research Letters.

\vspace{1ex}

\noindent [MP22b]\, S. Marouani, D. Popovici --- {\it Some Properties of Balanced Hyperbolic Compact Complex Manifolds} --- Internat. J. of Math., {\bf 33}, 3 (2022) 2250019, DOI : 10.1142/S0129167X22500197.

\vspace{1ex}

\noindent [MQ98]\, M. McQuillan --- {\it Diophantine Approximations and Foliations} --- Inst. Hautes Etudes Sci. Publ. Math., {\bf 87}, 121-174 (1998).

\vspace{1ex}

\noindent [MO22]\, S. Miebach, D.  Oeljeklaus --- {\it Compact Complex Non-K\"ahler Manifolds Associated With Totally Real Reciprocal Units} --- Math. Z. 301, No. 3, 2747-2760 (2022).

\vspace{1ex}

\noindent [Nak75]\, I. Nakamura --- {\it Complex Parallelisable Manifolds and Their Small Deformations} --- J. Differ. Geom. {\bf 10}, (1975), 85-112.

\vspace{1ex}

\noindent [OT05]\, K. Oeljekkaus, M. Toma --- {\it Non-K\"ahler Compact Complex Manifolds Associated to Number Fields} --- Ann. Inst. Fourier (Grenoble) {\bf 55} (2005), no. 1, 161-171. 

\vspace{1ex}

\noindent [Oti22]\, A Otiman --- {\it Special Hermitian Metrics on Oeljeklaus-Toma Manifolds} --- Bull. Lond. Math. Soc. {\bf 54} (2022), no. 2, 655-667.

\vspace{1ex}

\noindent [Pop15]\, D. Popovici --- {\it Aeppli Cohomology Classes Associated with Gauduchon Metrics on Compact Complex Manifolds} --- Bull. Soc. Math. France {\bf 143}, no. 4 (2015), p. 763-800.

\vspace{1ex}

\noindent [Pop23]\, D. Popovici --- {\it Pluriclosed Star Split Hermitian Metrics} ---  Math. Z. {\bf 305}, 7 (2023),

\noindent https://doi.org/10.1007/s00209-023-03344-0.

\vspace{1ex}

\noindent [TY09]\, N. Tsuchiya, A. Yamakawa --- {\it Lattices of Some Solvable Lie Groups and Actions of Products of Affine Groups} --- Tohoku Math. J. (2) {\bf 61} (2009), no. 3, 349-364.

\vspace{1ex}

\noindent [Tsu94]\, K. Tsukada --- {\it Holomorphic Forms and Holomorphic Vector Fields on Compact Generalized Hopf Manifolds} --- Compositio Math. {\bf 93} (1994), no. 1, 1-22.

\vspace{1ex}

\noindent [Vai79]\, I. Vaisman --- {\it Locally Conformal K\"ahler Manifolds With Parallel Lee Form} ---Rend. Mat. (6)  {\bf 12} (1979), no. 2, 263-284.

\vspace{1ex}

\noindent [Voi02]\, C. Voisin --- {\it Hodge Theory and Complex Algebraic Geometry. I.} --- Cambridge Studies in Advanced Mathematics, 76, Cambridge University Press, Cambridge, 2002.

%\noindent [Wan54]\, H.-C. Wang --- {\it Complex Parallisable Manifolds} --- Proc. Amer. Math. Soc. {\bf 5} (1954), 771--776.

%\vspace{1ex}

%\noindent [Yac98]\, A. Yachou --- {\it Sur les vari\'et\'es semi-k\"ahl\'eriennes} --- PhD Thesis, University of Lille.

\vspace{2ex}

\noindent Department of Mathematics    \hfill Institut de Math\'ematiques de Toulouse

\noindent Graduate School of Science    \hfill  Universit\'e Paul Sabatier

\noindent Osaka University, Osaka, Japan    \hfill  118 route de Narbonne, 31062 Toulouse, France

\noindent  Email: kasuya@math.sci.osaka-u.ac.jp                \hfill     Email: popovici@math.univ-toulouse.fr

\end{document}